\documentclass{amsart}

\usepackage{amssymb,amsmath}
\usepackage{comment}

\numberwithin{equation}{section}

\theoremstyle{plain}
\newtheorem{theorem}{Theorem}[section]
\newtheorem{proposition}[theorem]{Proposition}
\newtheorem{lemma}[theorem]{Lemma}
\newtheorem{corollary}[theorem]{Corollary}
\newtheorem{prob}[theorem]{Problem}
\newtheorem{alg}[theorem]{Algorithm}
\newtheorem{rmk}[theorem]{Remark}
\newtheorem*{claim*}{Claim}

\newenvironment{enumalph}
{\begin{enumerate}}
{\end{enumerate}}

\newenvironment{enumalgalph}
{\begin{enumerate}}
{\end{enumerate}}

\DeclareMathOperator{\Cl}{Cl}
\DeclareMathOperator{\disc}{disc}
\DeclareMathOperator{\End}{End}

\DeclareMathOperator{\N}{N}
\DeclareMathOperator{\nrd}{nrd}

\DeclareMathOperator{\rad}{rad}

\DeclareMathOperator{\SV}{SV}
\DeclareMathOperator{\trd}{trd}
\DeclareMathOperator{\Tr}{Tr}

\newcommand{\C}{\mathbb C}
\newcommand{\F}{\mathbb F}

\newcommand{\PP}{\mathbb P}
\newcommand{\Q}{\mathbb Q}
\newcommand{\R}{\mathbb R}
\newcommand{\Z}{\mathbb Z}

\newcommand{\calC}{\mathcal C}
\newcommand{\calE}{\mathcal E}
\newcommand{\calF}{\mathcal F}
\newcommand{\calI}{\mathcal I}
\newcommand{\calJ}{\mathcal J}
\newcommand{\calO}{\mathcal O}

\newcommand{\eps}{\epsilon}

\newcommand{\fraka}{\mathfrak{a}}
\newcommand{\frakb}{\mathfrak{b}}
\newcommand{\frakc}{\mathfrak{c}}
\newcommand{\frakd}{\mathfrak{d}}
\newcommand{\frakD}{\mathfrak{D}}
\newcommand{\frakf}{\mathfrak{f}}

\newcommand{\frakH}{\mathfrak{H}}
\newcommand{\frakl}{\mathfrak{l}}
\newcommand{\frakn}{\mathfrak{n}}
\newcommand{\frakN}{\mathfrak{N}}
\newcommand{\frakp}{\mathfrak{p}}
\newcommand{\frakq}{\mathfrak{q}}

\newcommand{\legen}[2]{\left(\frac{#1}{#2}\right)}

\newcommand{\quat}[2]{\displaystyle{\biggl(\frac{#1}{#2}\biggr)}}

\newcommand{\la}{\langle}
\newcommand{\ra}{\rangle}

\title[Quaternion ideal classes]{Algorithmic enumeration of ideal classes for quaternion orders}

\author{Markus Kirschmer}
\address{Lehrstuhl D f\"ur Mathematik, RWTH Aachen University, Templergraben 64, 52062 Aachen, Germany}
\email{Markus.Kirschmer@math.rwth-aachen.de}
\author{John Voight}
\address{Department of Mathematics and Statistics, University of Vermont, 16 Colchester Ave, Burlington, VT 05401, USA}
\email{jvoight@gmail.com}


\begin{document}

\begin{abstract}
We provide algorithms to count and enumerate representatives of the (right) ideal classes of an Eichler order in a quaternion algebra defined over a number field.  We analyze the run time of these algorithms and consider several related problems, including the computation of two-sided ideal classes, isomorphism classes of orders, connecting ideals for orders, and ideal principalization.  We conclude by giving the complete list of definite Eichler orders with class number at most $2$.
\end{abstract}

\maketitle

Since the very first calculations of Gauss for imaginary quadratic fields, the problem of computing the class group of a number field $F$ has seen broad interest.  Due to the evident close association between the class number and regulator (embodied in the Dirichlet class number formula), one often computes the class group and unit group in tandem as follows.

{\scshape Problem} (\textsf{ClassUnitGroup}($\Z_F$)). Given the ring of integers $\Z_F$ of a number field $F$, compute the class group $\Cl \Z_F$ and unit group $\Z_F^*$.

This problem appears in general to be quite difficult.  The best known (probabilistic) algorithm is due to Buchmann \cite{Buchmann}: for a field $F$ of degree $n$ and absolute discriminant $d_F$, it runs in time $d_F^{1/2}(\log d_F)^{O(n)}$ without any hypothesis \cite{LenstraRingInt}, and assuming the Generalized Riemann Hypothesis (GRH), it runs in \emph{expected} time $\exp\bigl(O\bigl((\log d_F)^{1/2}(\log\log d_F)^{1/2}\bigr)\bigr)$, where the implied $O$-constant depends on $n$. 

According to the Brauer-Siegel theorem, already the case of imaginary quadratic fields shows that the class group is often roughly as large as $d_F^{1/2}(\log d_F)^{O(1)}$.  Similarly, for the case of real quadratic fields, a fundamental unit is conjectured to have height often as large as $d_F^{1/2}(\log d_F)^{O(1)}$, so even to write down the output in a na\"ive way requires exponential time (but see Remark \ref{polysize}).  The problem of simply computing the class number $h(F)=\#\Cl \Z_F$, or for that matter determining whether or not a given ideal of $\Z_F$ is principal, appears in general to be no easier than solving Problem (\textsf{ClassUnitGroup}).

In this article, we consider a noncommutative generalization of the above problem.  We refer to \S 1 for precise definitions and specification of the input and output.

{\scshape Problem} (\textsf{ClassNumber}($\calO$)).
\emph{Given an Eichler order $\calO$ in a quaternion algebra over a number field $F$, compute the class number $h(\calO)$.}

{\scshape Problem} (\textsf{ClassSet}($\calO$)).
\emph{Given an Eichler order $\calO$ in a quaternion algebra over a number field $F$, compute a set of representatives for the set of invertible right $\calO$-ideal classes $\Cl \calO$.}

The main results of this article are embodied in the following two theorems, which provide algorithms to solve these two problems depending on whether the order is definite or indefinite.

\pagebreak[2]

{\scshape Theorem A}. 
\begin{enumalph}
\item \emph{If $\calO$ is indefinite, Problem \textup{(\textsf{ClassNumber}($\calO$))} is deterministic polynomial-time reducible to Problem \textup{(\textsf{ClassUnitGroup}($\Z_F$))}.}

\item \emph{If $\calO$ is definite, then Problem \textup{(\textsf{ClassNumber}($\calO$))} is reducible in probabilistic time 
\[ O\bigl(d_F^{3/2}\log^4 d_F+\log^2 \N \frakd\bigr) \] 
to the factorization of the discriminant $\frakd$ of $\calO$ and $O(2^n)$ instances of Problem \textup{\textsf{(ClassUnitGroup)}} with fields having discriminant of size $O(d_F^{5/2})$.}
\end{enumalph}

Here and throughout the paper, unless otherwise noted the implied $O$-constants are absolute.

{\scshape Corollary}.
\emph{There exists a probabilistic polynomial-time algorithm to solve Problem \textup{(\textsf{ClassNumber})} over a fixed field $F$ for indefinite orders and definite orders with factored discriminant.}

{\scshape Theorem B}.
\emph{There exists an algorithm to solve Problem \textup{(\textsf{ClassSet})} for orders over a fixed field $F$.  This algorithm runs in probabilistic polynomial time in the size of the output for indefinite orders and for definite orders with factored discriminant.}

It is important to note in Theorem B that we do not claim to be able to solve Problem (\textsf{ClassSet}) in probabilistic polynomial time in the size of the \emph{input}, since the output is of exponential size and so even to write ideal representatives (in the usual way) requires exponential time.  

The algorithmic results embodied in Theorems A and B have many applications.  Quaternion algebras are the noncommutative analogues of quadratic field extensions and so present an interesting and rewarding class to analyze.  For example, the norm form on a quaternion order gives rise to quadratic modules of rank $3$ and $4$ over $\Z_F$ and computing ideal classes amounts to finding all isometry classes of forms in the same genus (see e.g.\ Alsina-Bayer \cite[Chapter 3]{AlsinaBayer} for the case $F=\Q$).  Ideal classes in quaternion orders are also intimately related to automorphic forms.  In the simplest case where $F=\Q$, the ideal classes of a maximal order in a quaternion algebra of discriminant $p$ are in bijection with the set of supersingular elliptic curves in characteristic $p$.  This correspondence has been exploited by Pizer \cite{Pizer}, Kohel \cite{Kohel}, and others to explicitly compute spaces of modular forms over $\Q$.  By extension, one can compute with Hilbert modular forms over a totally real field $F$ of even degree using these methods via the Jacquet-Langlands correspondence \cite{Dembele}, and the algorithms described below have already been used for this purpose \cite{DembeleDonnelly}.  Finally, this work allows explicit computations with Shimura curves, including the computation of CM points \cite{VoightANTS}.

The outline of this article is as follows.  In Section 1, we review background material from the theory of Eichler orders of quaternion algebras and their ideals.  In Section 2, we introduce the algorithmic problems we will consider and discuss some of their interrelationships.  In Section 3, we treat the problem of computing the set of two-sided ideal classes and connecting ideals for Eichler orders.  In Section 4, we enumerate ideal classes in indefinite orders; we deduce Theorem A and its corollary in this case from Eichler's theorem of norms.  In Section 5, we introduce the Eichler mass formula which gives rise to an algorithm to count ideal classes in a definite quaternion order, completing the proof of Theorem A.  In Section 6, we discuss ideal principalization in definite orders, rigorously analyzing the lattice search employed to find a generator, if it exists.  In Section 7, we show how to enumerate ideal classes in definite orders, and use the Ramanujan property of the $\frakp$-neighbors graph to prove Theorem B.  Finally, in Section 8, we use our implementation of these algorithms in \textsf{Magma} \cite{Magma} to compute the complete list of definite quaternion orders (over an arbitrary totally real field $F$) with class number at most $2$ (Tables \ref{table1}--\ref{table2}): counted up to a natural notion of equivalence, there are exactly $74$ equivalence classes of definite Eichler orders with class number $1$ and $172$ with class number $2$.

We conclude this introduction by indicating two other results of significance found in the paper.  In Section 5, we give a rigorous complexity analysis for computing the value $\zeta_F(-1)$ of the Dedekind zeta function of a totally real field $F$.  We then
prove a complexity result (Proposition \ref{factorints}) which serves as a partial converse to Theorem B: the problem of factoring integers $a$ with $O(\log \log a)$ prime factors is probabilistic polynomial time reducible to Problem \textup{\textsf{(ClassNumber)}} over $\Q$.  In particular, if one can compute the class number of a definite rational quaternion order efficiently, then one can factor RSA moduli $a=pq$ efficiently.

The authors would like to thank the \textsf{Magma} group at the University of Sydney for their support, especially the assistance and advice of Steve Donnelly, as well as the reviewers for their very helpful comments and suggestions.  The authors would also like to thank Daniel Smertnig whose questions led to the discovery of an error in the proof of Theorem B (corrected here).  The second author would like to thank Reinier Br\"oker, Johannes Buchmann, Pete Clark, Henri Cohen, Tim Dokchitser, Claus Fieker, Eyal Goren, David Kohel, and Damien Stehl\'e for their valuable input.

\section{Quaternion algebras, orders, and ideals}

In this section, we introduce quaternion algebras and orders, and describe some of their basic properties; for further reading, see Reiner \cite{Reiner}, Vign\'eras \cite{Vigneras}, and Brezezi\'nski \cite{Brzezinski}.  Throughout, let $F$ be a number field of degree $[F:\Q]=n$ and absolute discriminant $d_F$, and let $\Z_F$ be its ring of integers.

\subsection*{Number rings}

We follow the usual algorithmic conventions for number fields and finitely generated abelian groups (see Cohen \cite{Cohen}, Lenstra \cite{LenstraRingInt}).  In particular, following Lenstra \cite[2.5]{LenstraRingInt}, to \emph{compute} a finitely generated abelian group $G$ means to specify a finite sequence $d_i \in \Z_{\geq 0}$ and an isomorphism $\bigoplus_i \Z/d_i \Z \xrightarrow{\sim} G$, in turn specified by the images of the standard generators.  Moreover, we represent a finitely generated torsion-free $\Z_F$-module $I$ by a \emph{pseudobasis} over $\Z_F$, writing $I = \bigoplus_i \fraka_i \gamma_i$ with $\fraka_i$ fractional ideals of $\Z_F$ and $\gamma_i \in I$.  See Cohen \cite[Chapter 1]{Cohen2} for methods of computing with finitely generated modules over Dedekind domains using pseudobases.

As in the introduction, we have the following basic problem.

\begin{prob}[\textsf{ClassUnitGroup}($\Z_F$)]
Given $\Z_F$, compute the class group $\Cl \Z_F$ and unit group $\Z_F^*$.
\end{prob}

\begin{rmk} \label{polysize}
{\rm The representation of the output of Problem \textup{(\textsf{ClassUnitGroup})} is not unique, and therefore different algorithms may produce correct output but conceivably of arbitrarily large size.  Indeed, we do not require that the outputted generators of the unit group $\Z_F^*$ to be represented in the usual way as a $\Z$-linear combination of an integral basis for $\Z_F$, since in general these elements can be of exponential size (as in the case of real quadratic fields).  Instead, we allow the units to be represented as a straight-line program involving elements of $\Z_F$ written in the usual way, for example as a (finite) product $\prod u^{c(u)}$ of elements $u \in F^*$ with $c(u) \in \Z$.  In this way, one may be able to write down a set of generators of subexponential size.}
\end{rmk}

\begin{proposition} \label{computeClZF}
There exists an algorithm to solve Problem \textup{(\textsf{ClassUnitGroup})} which runs in time $d_F^{1/2}(\log d_F)^{O(n)}$; assuming the generalized Riemann hypothesis \textup{(GRH)} and a ``smoothness condition'', this algorithm runs in time
\[ \exp\bigl(O((\log d_F)^{1/2}(\log\log d_F)^{1/2})\bigr), \] 
where the implied constant depends on $n$.
\end{proposition}

The algorithm underlying Proposition \ref{computeClZF} is due to Buchmann \cite{Buchmann} (see therein for the ``smoothness condition'', which is known to hold for quadratic fields).  See Lenstra \cite[Theorem 5.5]{LenstraRingInt}, Cohen-Diaz y Diaz-Olivier \cite{CDyDO}, and Cohen \cite[Algorithm 6.5.9]{Cohen} for further detail, and also Schoof \cite{Schoof} for a detailed analysis from the perspective of Arakelov geometry.

\begin{rmk}
{\rm A deterministic variant of the algorithm in Proposition \ref{computeClZF} runs in time $d_F^{3/4}(\log d_F)^{O(n)}$, due to the need to factor polynomials over finite fields.  We allow probabilistic algorithms in what follows.}
\end{rmk}

Further, there exists an algorithm which, given the internal calculations involved in the class group computation of Proposition \ref{computeClZF}, determines whether or not an ideal $\fraka \subset \Z_F$ is principal and, if so, outputs a generator (see Cohen \cite[Algorithm 6.5.10]{Cohen}).  No estimates on the running time of this algorithm have been proven, but it is reasonable to expect that they are no worse than the time for the class group computation itself.  (See also Remark \ref{princZFrmk} below for an alternative approach, which gives a principalization algorithm which runs in deterministic polynomial time over a fixed totally real field $F$.)

\subsection*{Quaternion algebras}

A \emph{quaternion algebra} $B$ over $F$ is a central simple algebra of dimension $4$ over $F$, or equivalently an $F$-algebra with generators $\alpha,\beta \in B$ such that 
\begin{equation} \label{quateq}
\alpha^2=a, \quad \beta^2=b, \quad \alpha\beta=-\beta\alpha
\end{equation}
with $a,b \in F^*$.  Such an algebra is denoted $B=\quat{a,b}{F}$ and is specified in bits by the elements $a,b \in F^*$, and an element $\gamma = x+y\alpha+z\beta+w\alpha\beta \in B$ is specified by the elements $x,y,z,w \in F$.

Let $B$ be a quaternion algebra over $F$.  Then $B$ has a unique (anti-)involution $\overline{\phantom{x}}:B \to B$ called \emph{conjugation} such that $\gamma+\overline{\gamma},\gamma\overline{\gamma} \in F$ for all $\gamma \in B$.  We define the \emph{reduced trace} and \emph{reduced norm} of $\gamma$ to be  $\trd(\gamma)=\gamma+\overline{\gamma}$ and $\nrd(\gamma)=\gamma\overline{\gamma}$, respectively.  For $B=\quat{a,b}{F}$ as in (\ref{quateq}) we have 
\begin{equation} \label{nrdtrd}
\overline{\gamma}=x-(y\alpha+z\beta+w\alpha\beta), \quad \trd(\gamma)=2x, \quad \nrd(\gamma)=x^2-ay^2-bz^2+abw^2.
\end{equation}

Let $K$ be a field containing $F$.  Then $B_K=B \otimes_F K$ is a quaternion algebra over $K$, and we say $K$ \emph{splits} $B$ if $B_K \cong M_2(K)$.  If $[K:F]=2$, then $K$ splits $B$ if and only if there exists an $F$-embedding $K \hookrightarrow B$.  

Now let $v$ be a place of $F$, and let $F_v$ denote the completion of $F$ at $v$.  We say $B$ is \emph{split} at $v$ if $F_v$ splits $B$, and otherwise we say that $B$ is \emph{ramified} at $v$.  The set of ramified places of $B$ is of even (finite) cardinality and uniquely characterizes $B$ up to $F$-algebra isomorphism.  We define the \emph{discriminant} $\frakD=\disc(B)$ of $B$ to be the ideal of $\Z_F$ given by the product of all finite ramified places of $B$.  One can compute the discriminant of $B$ in probabilistic polynomial time given an algorithm for integer factorization \cite{VoightHilbert}.

We say that $B$ is \emph{totally definite} if $F$ is totally real and every real infinite place of $F$ is ramified in $B$; otherwise, we say that $B$ is \emph{indefinite} or that $B$ \emph{satisfies the Eichler condition}.

\subsection*{Quaternion orders}

A $\Z_F$-\emph{lattice} $I \subset B$ is a finitely generated $\Z_F$-submodule of $B$ with $IF=B$.  An \emph{order} $\calO \subset B$ is a $\Z_F$-lattice which is also a subring of $B$ (hence $1 \in \calO$), and an order is \emph{maximal} if it is not properly contained in another order.  We represent a $\Z_F$-lattice by a pseudobasis, as above.  The problem of computing a maximal order $\calO$ is probabilistic polynomial-time equivalent to integer factorization \cite{VoightThesis}.  

We will use the general principle throughout that $\Z_F$-lattices are determined by their localizations.  For every prime ideal $\frakp$ of $\Z_F$, let $\Z_{F,\frakp}$ denote the completion of $\Z_F$ at $\frakp$.  For a $\Z_F$-lattice $I$, we abbreviate $I_\frakp = I \otimes_{\Z_F} \Z_{F,\frakp}$.  Then for two $\Z_F$-lattices $I,J \subset B$, we have $I=J$ if and only if $I_\frakp=J_\frakp$ for all primes $\frakp$.

For $\gamma_1,\dots,\gamma_4 \in B$, we let $\disc(\gamma_1,\dots,\gamma_4)=\det(\trd(\gamma_i\gamma_j))_{i,j=1,\dots,4}$.  
For an order $\calO \subset B$, the ideal of $\Z_F$ generated by the set 
\[ \{\disc(\gamma_1,\dots,\gamma_4) : \gamma_i \in \calO\} \]
is a square, and we define the \emph{(reduced) discriminant} $\frakd=\disc(\calO)$ of $\calO$ to be the square root of this ideal.  If $\calO=\bigoplus_i \fraka_i \gamma_i$ then we compute the discriminant as
\[ \frakd^2=\left(\fraka_1 \cdots \fraka_4\right)^2 \disc(\gamma_1,\dots,\gamma_4). \]
An order $\calO$ is maximal if and only if $\frakd=\frakD$.  

An \emph{Eichler order} is the intersection of two maximal orders, and it is this class of orders which we will study throughout.  The \emph{level} of an Eichler order $\calO$ is the ideal $\frakN \subset \Z_F$ satisfying $\frakd=\frakD\frakN$; the level $\frakN$ is coprime to the discriminant $\frakD$ of $B$.  Alternatively, given a maximal order $\calO \subset B$, an ideal $\frakN$ coprime to $\frakD$ and an embedding $\iota_\frakN:\calO \hookrightarrow M_2(\Z_{F,\frakN})$ where $\Z_{F,\frakN}$ denotes the completion of $\Z_F$ at $\frakN$, an Eichler order of level $\frakN$ is given by
\begin{equation} \label{O0(N)}
\calO_0(\frakN)=\left\{\gamma \in \calO : \iota_\frakN(\gamma) \text{ is upper triangular modulo $\frakN$}\right\},
\end{equation}
and all Eichler orders arise in this way up to conjugation.  In particular \cite[Theorem 39.14]{Reiner}, an order $\calO$ is hereditary (all one-sided ideals of $\calO$ are projective) if and only if $\calO$ is an Eichler order with squarefree level.  

We can compute an explicit pseudobasis for an Eichler order $\calO_0(\frakN)$ from the description (\ref{O0(N)}) as follows.  First, we compute a maximal order $\calO \subset B$ as above.  Next, for each prime power $\frakp^e \parallel \frakN$, we compute an embedding $\iota_\frakp:\calO \hookrightarrow M_2(\Z_{F,\frakp})$; this can be accomplished in probabilistic polynomial time \cite{VoightHilbert}.  From $\iota_\frakp$, one easily computes $\calO_0(\frakp^e)$ using linear algebra.  Then $\calO_0(\frakN) = \bigcap_{\frakp^e \parallel N} \calO_0(\frakp^e)$, and this intersection can be computed as $\Z_F$-lattices.

Two orders $\calO,\calO'$ are \emph{conjugate} (also \emph{isomorphic} or \emph{of the same type}) if there exists $\nu \in B^*$ such that $\calO'=\nu^{-1}\calO\nu$, and we write $\calO \cong \calO'$.

\begin{proposition}[{\cite[Corollaire III.5.5]{Vigneras}}]
The number of isomorphism classes of Eichler orders $\calO \subset B$ of level $\frakN$ is finite.  
\end{proposition}

\subsection*{Quaternion ideals}

We define the \emph{reduced norm} $\nrd(I)$ of a $\Z_F$-lattice $I$ to be the fractional ideal of $\Z_F$ generated by the set $\{\nrd(\gamma) : \gamma \in I\}$.  

Let $I,J$ be $\Z_F$-lattices in $B$.  We define the product $IJ$ to be the $\Z_F$-submodule of $B$ generated by the set $\{\alpha\beta : \alpha \in I,\ \beta \in J\}$; we have $\nrd(IJ) \supset \nrd(I)\nrd(J)$.  We define the \emph{left colon} 
\[ (I:J)_L = \{\gamma \in B : \gamma J \subset I\} \]
and similarly the \emph{right colon}
\[ (I:J)_R = \{\gamma \in B : J \gamma \subset I\}. \]
The colons $(I:J)_L,(I:J)_R$ and the product $IJ$ are $\Z_F$-lattices.  If $I=J$, then $(I:I)_L=\calO_L(I)$ (resp.\ $(I:I)_R=\calO_R(I)$) also has the structure of a ring, called the \emph{left} (resp.\ \emph{right}) \emph{order} of the $\Z_F$-lattice $I$.  One can compute the left and right colon in deterministic polynomial time using the Hermite normal form for $\Z_F$-lattices (see Friedrichs \cite[\S 2.3]{Friedrichs}).  

Let $\calO \subset B$ be an order.  A \emph{right fractional $\calO$-ideal} is a $\Z_F$-lattice $I$ such that $\calO_R(I)=\calO$.  In a similar fashion, we may define left fractional ideals; however, conjugation 
\[ I \mapsto \overline{I} = \{\overline{\gamma} : \gamma \in I\} \]
gives a bijection between the sets of right and left fractional $\calO$-ideals, so when dealing with one-sided fractional ideals it suffices to work with right fractional ideals.  If $I$ (resp.\ $J$) is a right fractional $\calO$-ideal then $\calO_R((I:J)_R) \subset \calO$ (resp.\ $\calO_L((I:J)_R) \subset \calO$).  Note that any $\Z_F$-lattice $I$ is by definition a right fractional $\calO_R(I)$-ideal (and left fractional $\calO_L(I)$-ideal).

A $\Z_F$-lattice $I$ is \emph{integral} if $I \subset \calO_R(I)$, or equivalently if $I$ is a right ideal of $\calO_R(I)$ in the usual sense; for any $\Z_F$-lattice $I$, there exists a nonzero $d \in \Z_F$ such that $dI$ is integral.

A $\Z_F$-lattice $I$ is a left fractional $\calO_L(I)$-ideal and a right fractional $\calO_R(I)$-ideal, and we say that $I$ is a fractional $\calO_L(I),\calO_R(I)$-ideal; if $\calO_L(I)=\calO_R(I)=\calO$ we say that $I$ is a \emph{two-sided} $\calO$-ideal.  

A right fractional $\calO$-ideal is \emph{left invertible} if there exists a left fractional $\calO$-ideal $I^{-1}$ such that $I^{-1}I=\calO$.  If $I$ is left invertible, then necessarily 
\[ I^{-1}=(\calO:I)_L=\overline{I}/\nrd(I). \]  
Equivalently, $I$ is left invertible if and only if $I$ is \emph{locally principal}, i.e., for each (finite) prime ideal $\frakp$ of $\Z_F$, the ideal $I_\frakp$ is a principal right $\calO_\frakp$-ideal.  It follows that if $I$ is an $\calO',\calO$-ideal then $I$ is left invertible if and only if $I$ is right invertible (analogously defined), and so we simply say $I$ is \emph{invertible}, and then $II^{-1}=\calO'$ and $I^{-1}=(\calO':I)_R$.  If $I$ is an invertible right fractional $\calO$-ideal and $J$ is an invertible left fractional $\calO$-ideal then $(IJ)^{-1}=J^{-1} I^{-1}$ and $\nrd(IJ)=\nrd(I)\nrd(J)$, and moreover
\[ \calO_L(IJ) = \calO_L(I) \text{ and } \calO_R(IJ) = \calO_R(J). \]  
We note that for an order $\calO$, every right fractional $\calO$-ideal $I$ is invertible if and only if $\calO$ is hereditary.  

Let $I,J$ be invertible right fractional ideals.  Then $(I:J)_R$ is a fractional $\calO_R(J),\calO_R(I)$-ideal and similarly $(I:J)_L$ is a fractional $\calO_L(I), \calO_L(J)$-ideal, and  so we will also call $(I:J)_L$ (resp.\ $(I:J)_R$) the \emph{left} (resp.\ \emph{right}) \emph{colon fractional ideal}.  

Let $I,J$ be invertible right fractional $\calO$-ideals.  We say that $I$ and $J$ are in the same \emph{right ideal class}, and write $I \sim J$, if there exists an $\alpha \in B^*$ such that $I=\alpha J$.  We have $I \sim J$ if and only if $I$ and $J$ are isomorphic as right $\calO$-modules, and so in this case we also say that $I$ and $J$ are \emph{isomorphic}.  It is clear that $\sim$ defines an equivalence relation on the set of right fractional ideals of $\calO$; we write $[I]$ for the ideal class of $I$.  Since $B$ is noncommutative, the ideal class $[IJ]$ of two right fractional $\calO$-ideals $I,J$ is in general not determined by the ideal classes $[I]$ and $[J]$, so the set of right ideal classes may not form a group.  We denote the set of invertible right $\calO$-ideal classes by $\Cl \calO$.  

The set of invertible two-sided fractional $\calO$-ideals forms a group under multiplication, and the quotient of this group by the (normal) subgroup of principal two-sided fractional $\calO$-ideals is called the \emph{two-sided ideal class group} of $\calO$; two invertible two-sided fractional $\calO$-ideals $I,J$ are in the same \emph{ideal class} if they are equal in the two-sided ideal class group of $\calO$, or equivalently if $IJ^{-1}$ is a principal two-sided fractional $\calO$-ideal.

An order $\calO$ is \emph{connected} to an order $\calO'$ if there exists an invertible fractional $\calO,\calO'$-ideal $I$, the \emph{connecting ideal}.  The relation of being connected is an equivalence relation, and two Eichler orders $\calO,\calO'$ are connected if and only if they have the same level $\frakN$.

\begin{proposition}[{\cite[Th\'eor\`eme III.5.4]{Vigneras}}, {\cite[\S 26]{Reiner}}]
The set $\Cl \calO$ is finite and $\#\Cl \calO$ is independent of the choice of Eichler order $\calO$ of a given level.
\end{proposition}

We let $h(\calO)=\#\Cl \calO$ denote the \emph{class number} of the Eichler order $\calO$.  

\section{Algorithmic problems}

In the remainder of this article, we will be concerned with a constellation of interrelated algorithmic problems which we now introduce.

\begin{prob}[\textsf{ClassNumber}($\calO$)] \label{computeclassnumber}
Given an Eichler order $\calO$, compute the class number $h(\calO)$.
\end{prob}

\begin{prob}[\textsf{ClassSet}($\calO$)] \label{computeidealclasses}
Given an Eichler order $\calO$, compute a set of representatives for the set of invertible right $\calO$-ideal classes $\Cl \calO$.
\end{prob}

Obviously, a solution to Problem \ref{computeidealclasses} (\textsf{ClassSet}) gives a solution to Problem \ref{computeclassnumber} (\textsf{ClassNumber}), but as we will see, this reduction is not the most efficient approach.

Given a set of representatives for $\Cl \calO$ and a right fractional ideal $I$ of $\calO$, we may also like to determine its class $[I] \in \Cl \calO$ and so we are led to the following problems.

\begin{prob}[\textsf{IsIsomorphic}($I,J$)] \label{computeisomorphism}
Given two invertible right fractional ideals $I,J$ of an Eichler order $\calO$, determine if $I \sim J$; and, if so, compute $\xi \in B^*$ such that $I=\xi J$.  
\end{prob}

\begin{prob}[\textsf{IsPrincipal}($I$)] \label{computeprincipal}
Given an invertible right fractional ideal $I$ of an Eichler order $\calO$, determine if $I$ is principal; and, if so, compute a generator $\xi$ of $I$.  
\end{prob}

In fact, these two problems are computationally equivalent.

\begin{lemma} \label{reductionisomtest}
Problem \textup{(\textsf{IsIsomorphic})} is equivalent to Problem \textup{(\textsf{IsPrincipal})}.
\end{lemma}

\begin{proof}
Let $I,J$ be invertible right fractional $\calO$-ideals.  Then $I=\xi J$ for $\xi \in B^*$ if and only if the left colon ideal $(I:J)_L$ is generated by $\xi$ as a right fractional $\calO_L(J)$-ideal.  Therefore, $I \sim J$ if and only if $(I:J)_L$ is a principal fractional right $\calO_L(J)$-ideal.
\end{proof}

We also have the corresponding problem for two-sided ideals.

\begin{prob}[\textsf{TwoSidedClassSet}($\calO$)] \label{twosidedclasses}
Given an Eichler order $\calO$, compute a set of representatives for the two-sided invertible ideal classes of $\calO$.
\end{prob}

Finally, we consider algorithmic problems for orders.

\begin{prob}[\textsf{IsConjugate}($\calO,\calO'$)] \label{conjugacytest}
Given two Eichler orders $\calO,\calO'$ of $B$, determine if $\calO \cong \calO'$; and, if so, compute $\nu \in B^*$ such that $\nu \calO \nu^{-1} = \calO'$.
\end{prob}

\begin{prob}[\textsf{ConjClassSet}($\calO$)] \label{conjugacyclasses}
Given an Eichler order $\calO$ of level $\frakN$, compute a set of representatives for the conjugacy classes of Eichler orders of level $\frakN$.
\end{prob}

\begin{prob}[\textsf{ConnectingIdeal}($\calO,\calO'$)] \label{connectingideal}
Given Eichler orders $\calO,\calO'$, compute a connecting ideal $I$ with $\calO_R(I)=\calO$ and $\calO_L(I)=\calO'$.  
\end{prob}

We conclude by relating Problem \ref{computeidealclasses} (\textsf{ClassSet}) to Problem \ref{conjugacyclasses} (\textsf{ConjClassSet}).  

\begin{proposition} \label{conjareequiv}
Let $\calO_i$ be representatives of the isomorphism classes of Eichler orders of level $\frakN$.  For each $i$, let $I_i$ be a connecting fractional $\calO_i,\calO$-ideal and let $J_{i,j}$ be representatives of the two-sided invertible fractional $\calO_i$-ideal classes.  

Then the set $\{J_{i,j} I_i\}_{i,j}$ is a complete set of representatives of $\Cl \calO$.
\end{proposition}

\begin{proof}
Let $I$ be an invertible right fractional $\calO$-ideal.  Then $\calO_L(I) \cong \calO_i$ for a uniquely determined $i$, so $\calO_L(I)=\nu^{-1} \calO_i \nu$ for some $\nu \in B^*$.  But then $I_i = \nu K I$ where $K=\nu^{-1} I_i I^{-1}$ is a two-sided invertible fractional $\calO_L(I)$-ideal, and so $I \sim K I_i \sim J_{i,j} I_i$ for some $j$, again uniquely determined.
\end{proof}

It follows from Proposition \ref{conjareequiv} that if one can solve Problem \ref{conjugacyclasses} (\textsf{ConjClassSet}) then one can solve Problem \ref{computeidealclasses} (\textsf{ClassSet}), given algorithms to solve Problem \ref{connectingideal} (\textsf{ConnectingIdeal}) and Problem \ref{twosidedclasses} (\textsf{TwoSidedClassSet}).  We will discuss this further in \S \S 3--4.

Conversely, if one can solve Problem \ref{computeidealclasses} (\textsf{ClassSet}) then one can solve Problem \ref{conjugacyclasses} (\textsf{ConjClassSet}) given an algorithm to solve Problem \ref{conjugacytest} (\textsf{IsConjugate}): indeed, by Proposition \ref{conjareequiv}, one obtains a set of representatives for the conjugacy classes of orders by computing $\calO_L(I)$ for $[I] \in \Cl \calO$.  The difficulty of solving Problem \ref{conjugacytest} (\textsf{IsConjugate}) is discussed in \S 4 and \S 6.

\section{Two-sided ideal classes and connecting ideals}

In this section, we discuss Problem \ref{twosidedclasses} (\textsf{TwoSidedClassSet}) and Problem \ref{connectingideal} (\textsf{ConnectingIdeal}).

\subsection*{Two-sided ideal classes}

Let $\calO \subset B$ be an Eichler order of discriminant $\frakd$ and level $\frakN$.  The two-sided ideals of $\calO$ admit a local description, as follows.  Let $F_\frakp$ denote the completion of $F$ at $\frakp$, let $\Z_{F,\frakp}$ denote its ring of integers, and let $\pi$ be a uniformizer for $\Z_{F,\frakp}$.  

First, suppose that $B_\frakp = B \otimes_F F_\frakp$ is a division ring.  Then the discrete valuation $v$ of $\Z_{F,\frakp}$ extends to $B_\frakp$, and $\calO_\frakp$ is the unique maximal order of $B_\frakp$.  The fractional right $\calO_\frakp$-ideals form a cyclic group generated by the principal ideal 
\[ \rad(\calO_\frakp)=\{\gamma \in \calO_\frakp : v(\gamma)>0\}; \]
in particular, they are all two-sided \cite[Theorem 13.2]{Reiner} and invertible.  We have
$\rad(\calO_\frakp)=[\calO_\frakp,\calO_\frakp]$, where $[\calO_\frakp,\calO_\frakp]$ denotes the two-sided $\calO_\frakp$-ideal generated by the set $\{\gamma\delta-\delta\gamma : \gamma,\delta \in \calO_\frakp \}$ (see Reiner \cite[Theorem 14.5]{Reiner}).  

Next, suppose that $B_\frakp \cong M_2(F_\frakp)$ and that
\begin{equation} \label{standardeichler}
\calO_\frakp \cong \begin{pmatrix} \Z_{F,\frakp} & \Z_{F,\frakp} \\ \frakp^e\Z_{F,\frakp} & \Z_{F,\frakp} \end{pmatrix} \subset M_2(\Z_{F,\frakp}),
\end{equation}
so that $\frakp^e \parallel \frakN$.  Then the principal (equivalently, invertible) two-sided fractional ideals of $\calO_\frakp$ form an abelian group generated by
$\pi\calO_\frakp$ and 
\[ \left(\begin{matrix} 0 & 1 \\ \pi^e & 0 \end{matrix}\right) \calO_\frakp = [\calO_\frakp,\calO_\frakp] \]
(see the proof given by Eichler \cite[Satz 5]{Eichler} for $e=1$, which applies \emph{mutatis mutandis} for all $e$). Since $[\calO_\frakp,\calO_\frakp]^2 = \pi^e \calO_\frakp$, this group is cyclic if and only if $e$ is odd or $e=0$.

In particular, it follows from the preceding discussion that $[\calO,\calO]$ is an invertible two-sided $\calO$-ideal, and we  have the following description of the group of two-sided ideals.

\begin{lemma}\label{factor_twosided}
The set of invertible fractional two-sided $\calO$-ideals forms an abelian group generated by the set
\[ \{\frakp \calO : \frakp \subset \Z_F \} \cup \{ [\calO,\calO] + \frakp^e\calO : \frakp^e \parallel \frakd \}. \]
\end{lemma}

\begin{proof} 
The group of invertible two-sided fractional $\calO$-ideals is abelian since it is so locally by the above. 

Let $I$ be an invertible two-sided fractional $\calO$-ideal.  Clearing denominators, we may assume $I$ is integral.  Let $M$ be an invertible maximal two-sided $\calO$-ideal containing $I$.  Then by maximality, there exists a unique prime ideal $\frakp$ of $\Z_F$ such that $M_\frakq = \calO_\frakq$ for all $\frakq \neq \frakp$.  Thus by the preceding discussion, $M = \frakp \calO$ or $M = [\calO,\calO] + \frakp^e \calO$ with $\frakp^e \parallel \frakd$.  Now $IM^{-1}$ is again integral and $\nrd(IM^{-1})=\nrd(I)/\nrd(M) \mid \nrd(I)$, so the result follows by induction.
\end{proof}

For an alternative proof of Lemma \ref{factor_twosided} when $\calO$ is hereditary, see Vigneras \cite[Th\'eor\`eme I.4.5]{Vigneras}.

\begin{proposition} \label{alltwosidedclasses}
The group of invertible, two-sided fractional ideal classes of $\calO$ is a finite abelian group generated by the classes of
\[ \{ \fraka \calO : [\fraka] \in \Cl(\Z_F) \} \cup \{ [\calO,\calO] + \frakp^e\calO : \frakp^e \parallel \frakd \}. \]
If $B$ is indefinite, one can omit all generators $[\calO,\calO] + \frakp^e\calO$ for which $e$ is even.
\end{proposition}

\begin{proof} 
Since the principal two-sided fractional $\calO$-ideals form a subgroup, the first statement follows from the preceding lemma. For the second statement, we skip ahead and apply Proposition \ref{ECcond}: if $e$ is even, the ideals $\frakp^{e/2}\calO$ and $[\calO,\calO] + \frakp^{e}\calO$ have the same reduced norm $\frakp^e$, so they are in the same ideal class if $B$ is indefinite.
\end{proof}

\begin{corollary} \label{twosidedsolve}
Problem \textup{(\textsf{TwoSidedClassSet}($\calO$))} for an Eichler order $\calO$ with factored discriminant $\frakd$ is polynomial-time reducible to Problem \textup{(\textsf{ClassUnitGroup})($\Z_F$)} and $O(h(\Z_F)^2 \N\frakd^2)$ instances of Problem \textup{(\textsf{IsIsomorphic})}.
\end{corollary}

\begin{proof}
Since $([\calO,\calO]+\frakp^e\calO)^2 = \frakp^e \calO$ (as this is true locally), there are at most $2^{\omega(\frakd)}h(\Z_F)$ two-sided ideal classes, where $\omega(\frakd)$ denotes the number of prime factors of $\frakd$.  We have trivially $2^{\omega(\frakd)} \leq \N \frakd$, and the result follows.
\end{proof}

\subsection*{Eichler orders and connecting ideals}

We now exhibit an algorithm to test if an order is an Eichler order.  

\begin{alg} \label{iseichler}
Let $\calO \subset B$ be an order of discriminant $\frakd=\frakD\frakN$ with $\frakN$ prime to the discriminant $\frakD$ of $B$, and let $\iota_\frakN:\calO \hookrightarrow M_2(F_{\frakN})$.  This algorithm determines if $\calO$ is an Eichler order, and if so, returns an element $\nu \in B$ such that $\iota_\frakN(\nu^{-1}\calO\nu)=\calO_0(\frakN)$ (as in (\ref{O0(N)})).
{\rm \begin{enumerate}
\item Compute $\mu \in B$ such that $\iota_\frakN(\mu^{-1}\calO\mu) \subset M_2(\Z_{F_\frakN})$.  Let $\iota_\frakN'=\mu \iota_\frakN \mu^{-1}$.
\item Factor the ideal $\frakN$, and for each prime power $\frakp^e \parallel \frakN$:
\begin{enumalgalph}
\item From the restriction $\iota_\frakp':\calO \hookrightarrow M_2(F_\frakp)$ of $\iota_\frakN'$, use linear algebra over $\Z_{F,\frakp}$ to test if there is a common eigenvector $(x_\frakp,z_\frakp) \in (\Z_F/\frakp^e)^2$ for the elements of a $\Z_{F,\frakp}$-module basis of $\calO_\frakp$.  If not, return \textup{\textsf{false}}.
\item Compute $y_\frakp,w_\frakp$ such that $N_\frakp=\begin{pmatrix} x_\frakp & y_\frakp \\ z_\frakp & w_\frakp \end{pmatrix} \in GL_2(\Z_F/\frakp^e)$.  
\end{enumalgalph}
\item By the Chinese remainder theorem, compute $\nu \in B$ such that $\iota_\frakp(\nu_\frakp) \equiv N_\frakp \pmod{\frakp^e}$.  Return \textup{\textsf{true}} and $\mu\nu$.
\end{enumerate}}
\end{alg}

\begin{proof}[Proof of correctness]
We refer to work of the second author \cite{VoightHilbert} for more on Step 1.

For the rest of the algorithm, we note that the property of being an Eichler order is local: in particular, we see that a local order $\calO_\frakp=\calO \otimes_{\Z_F} \Z_{F,\frakp}$ with $\disc(\calO_\frakp)=\frakp^e$ is Eichler if and only if there exists such a common eigenvector modulo $\frakp^e$ of all $\gamma \in \calO_\frakp$.  Conjugation by the matrix $N_\frakp$ as in Step 2b then shows that $(1,0)$ is an eigenvector modulo all such $\frakp^e$, as desired.
\end{proof}

Now let $\calO, \calO'$ be two Eichler orders in $B$ having the same level $\frakN$.  We consider Problem \ref{connectingideal} (\textsf{ConnectingIdeal}) and compute an invertible $\calO',\calO$-ideal $I$.

For any prime $\frakp \nmid \frakN$, by maximality the $\Z_{F,\frakp}$-lattice $(\calO'\calO)_\frakp$ is a $\calO_\frakp',\calO_\frakp$-ideal.  So suppose $\frakp \mid \frakN$.  Since any two Eichler orders of the same level are locally conjugate, there exists $\nu_\frakp \in B_\frakp$ such that $\calO_\frakp'=\nu_\frakp \calO_\frakp \nu_\frakp^{-1}$.  The map $I \mapsto \nu_\frakp I$ gives an equivalence between the category of fractional two-sided $\calO_\frakp$-ideals and the category of fractional $\calO_\frakp',\calO_\frakp$-ideals.  

From this analysis, we arrive at the following algorithm.

\begin{alg} \label{connidealalg}
Let $\calO,\calO' \subset B$ be Eichler orders of level $\frakN$.  This algorithm computes an invertible fractional $\calO',\calO$-ideal.
{\rm \begin{enumerate}
\item Let $\nu,\nu'$ be the output of Algorithm \ref{iseichler} for the orders $\calO,\calO'$ and a common choice of splitting $\iota_\frakN$. 
\item Compute a nonzero $d \in \Z_F$ such that $\mu := d\nu^{-1}\nu' \in \calO'$ as follows: write $\nu^{-1}\nu'$ in terms of a $\Z_F$-pseudobasis for $\calO'$, and compute a nonzero $d$ as the least common multiple of the denominators of the coefficients of $\nu^{-1}\nu'$.  
\item Compute $\nrd(\mu)\Z_F=\frakn \fraka$ with $\fraka$ prime to $\frakN$, and return the $\Z_F$-lattice $I := \mu\calO + \frakn \calO'\calO$.
\end{enumerate}}
\end{alg}

\begin{proof}[Proof of correctness]
In Step 1, we obtain from Algorithm \ref{iseichler} that for all $\frakp \mid \frakN$ we have $\calO_\frakp=\nu_\frakp\calO_0(\frakN)_\frakp\nu_\frakp^{-1}$ and $\calO_\frakp=\nu_\frakp'\calO_0(\frakN)_\frakp\nu_\frakp'^{-1}$.  It is clear that Step 2 gives the correct output, and hence $\calO'_\frakp=\mu_\frakp\calO_\frakp\mu_\frakp^{-1}$ for all such $\frakp$.

To conclude, we need to show that $\calO_L(I)=\calO'$ and $\calO_R(I)=\calO$.  It suffices to check this locally.  For any prime $\frakp \nmid \frakN$, we have $\mu_\frakp \in \calO_\frakp'=\frakn_\frakp\calO'_\frakp$, so $I_\frakp=(\calO'\calO)_\frakp$, which is a fractional $\calO_\frakp',\calO_\frakp$-ideal by the above.  For $\frakp \mid \frakN$, we have 
\[ \frakn_\frakp \calO_\frakp' = \calO_\frakp' \frakn_\frakp = \calO'_\frakp \overline{\mu}_\frakp \mu_\frakp \subset \calO'_\frakp \mu_\frakp = \mu_\frakp \calO_\frakp \]
since $\overline{\mu}_\frakp \in \calO'_\frakp$.  Hence $I_\frakp=\mu_\frakp \calO_\frakp$ and the result follows by the equivalence above, since $\calO_\frakp$ is a two-sided $\calO_\frakp$-ideal so $I_\frakp$ is a fractional $\calO_\frakp',\calO_\frakp$-ideal.
\end{proof}

\begin{corollary} \label{isconjbyconnectingideal}
Problem \textup{(\textsf{IsConjugate})} for two Eichler orders with factored discriminant $\frakd$ is probabilistic polynomial-time reducible to Problems \textup{(\textsf{TwoSidedClassSet})} and \textup{(\textsf{IsIsomorphic})}.
\end{corollary}

\begin{proof}
% By Proposition \ref{conjareequiv}, any two connecting ideals for $\calO,\calO'$ differ by a two-sided ideal, and in particular the orders $\calO,\calO'$ are isomorphic if and only if a connecting ideal is isomorphic to a two-sided ideal.  
By Proposition \ref{conjareequiv}, if $I$ is an invertible right $\calO$-ideal, then $\calO_L(I)$ is conjugate to $\calO$ if and only if $I$ is equivalent to an invertible two-sided $\calO$-ideal $J$: in fact, if $J=\nu I$ then $\calO = \nu\calO'\nu^{-1}$.  Thus, to check whether two given Eichler orders $\calO, \calO'$ (of the same level) are conjugate, it suffices to construct a connecting ideal $I$ as in Algorithm \ref{connidealalg}---which can be done in probabilistic polynomial time---and one can accordingly check for an isomorphism given a solution to Problem (\textsf{TwoSidedClassSet}) and (\textsf{IsIsomorphic}).
\end{proof}

\section{Ideal classes in indefinite orders}

In this section, we discuss the difficulty of solving Problems \ref{computeclassnumber} (\textsf{ClassNumber}) and \ref{computeidealclasses} (\textsf{ClassSet}) in the indefinite case.  

Let $B$ be an indefinite quaternion algebra and let $\calO \subset B$ be an Eichler order.  Let $S_\infty$ denote the set of ramified (real) infinite places of $B$ and let $\Cl_{S_\infty} \Z_F$ denote the ray class group of $\Z_F$ with modulus $S_\infty$.  The quotient group $\Cl_{S_\infty} \Z_F/\Cl \Z_F$ is an elementary $2$-group isomorphic to $\Z_{F,S_\infty}^*/\Z_F^{*2}$ where 
\[ \Z_{F,S_\infty}^*=\{u \in \Z_F^*:v(u)>0 \text{ for all $v \in S_\infty$}\}.\]

\begin{proposition}[Eichler's theorem] \label{ECcond}
If $B$ is indefinite and $\calO \subset B$ is an Eichler order, then the reduced norm map
\[ \nrd:\Cl \calO \to \Cl_{S_\infty} \Z_F \]
is a bijection (of sets).
\end{proposition}

For a proof of this proposition, see Reiner {\cite[Corollary 34.21]{Reiner} or Vign\'eras \cite[Th\'eor\`eme III.5.7]{Vigneras}.  We have the following immediate corollary, which proves Theorem A in the indefinite case; we restate it here for convenience.

% \begin{rmk}
% If $B$ is a semisimple algebra over $F$ which satisfies the Eichler condition, i.e., $B$ has no simple component which is a definite quaternion algebra, then in fact Proposition \ref{ECcond} holds for any maximal order $\calO \subset B$.
% \end{rmk}

\begin{corollary} \label{indefO(1)}
If $B$ is indefinite, then Problem \textup{(\textsf{ClassNumber}($\calO$))} is reducible in deterministic polynomial time to Problem \textup{(\textsf{ClassUnitGroup}($\Z_F$))}.  
\end{corollary}

In other words, there exists an algorithm to solve Problem (\textsf{ClassNumber})} which, given an algorithm to solve Problem \textup{\textsf{(ClassUnitGroup)}}, runs in deterministic polynomial time in its output size.  (See Remark \ref{polysize}.)

\begin{proof}
We compute $h(\calO)=\#(\Z_{F,S_\infty}^*/\Z_F^{*2}) h(\Z_F) $.  Given the class group $\Cl \Z_F$, its order $h(\Z_F)=\#\Cl \Z_F$ can be computed in polynomial time.  Similarly, given generators for the unit group $\Z_F^*$, one can compute in deterministic polynomial time (in the size of their representation) their signs for each real place $v$, and using linear algebra over $\F_2$ determine the $2$-rank of the group $\Z_{F,S_\infty}^*/\Z_F^{*2}$.
\end{proof}

It follows immediately from Corollary \ref{indefO(1)} that for Eichler orders over a fixed number field $F$, Problem \textup{\ref{computeclassnumber} (\textsf{ClassNumber})} can be solved in time $O(1)$, which proves the corollary to Theorem A.  

Next, we discuss Theorem B in the indefinite case.  First, we exhibit an auxiliary algorithm for computing ideals with given norm, which works for both definite and indefinite quaternion orders.

\begin{alg} \label{IdealOfNorm}
Let $\calO$ be an Eichler order and let $\fraka \subset \Z_F$ be an ideal which is coprime to $\frakD$.  This algorithm returns an invertible right $\calO$-ideal $I$ such that $\nrd(I)=\fraka$.
{\rm \begin{enumerate}
\item Factor $\fraka$ into prime ideals.
\item For each $\frakp^e \parallel \fraka$, find a zero of the quadratic form $\nrd(\alpha_\frakp) \equiv 0 \pmod{\frakp}$, and choose a random lift of $\alpha_\frakp$ modulo $\frakp^2$ so that $\nrd(\alpha_\frakp)$ is a uniformizer at $\frakp$.  Let $\beta_\frakp=\alpha_\frakp^e$.
\item Use the Chinese remainder theorem to find $\beta \in \Z_F$ such that $\beta \equiv \beta_\frakp \pmod{\frakp^e}$ for all $\frakp^e \parallel \fraka$.  Return the right $\calO$-ideal $\beta\calO + \fraka\calO$.
\end{enumerate}}
\end{alg}

From Eichler's theorem (Proposition \ref{ECcond}), we then have the following straightforward algorithm.

\begin{alg} \label{ClassSetEichler}
Let $\calO$ be an indefinite Eichler order.  This algorithm solves Problem \ref{computeidealclasses} \textup{(\textsf{ClassSet})}.
{\rm \begin{enumerate}
\item For each $\fraka$ in a set of representatives for $\Cl_{S_\infty} \Z_F/2\Cl_{S_\infty} \Z_F$, compute using Algorithm \ref{IdealOfNorm} an ideal $I_\fraka$ of norm $\fraka$.
\item Return the set $\{\frakc I_\fraka\}_{\fraka,\frakc}$, with $\frakc^2$ in a set of representatives of $2\Cl_{S_\infty} \Z_F$.
\end{enumerate}}
\end{alg}

\begin{proposition} \label{indefclassset}
Problem \textup{\ref{computeidealclasses} (\textsf{ClassSet})} for indefinite orders over a fixed field $F$ can be solved in probabilistic polynomial time.  
\end{proposition}

\begin{proof}
One can solve Problem \textup{\textsf{(ClassUnitGroup)}} for the field $F$ in constant time and one can further factor the generating ideals $\fraka$ which are given as output.  The statement follows by noting that Step 2 of Algorithm \ref{IdealOfNorm} can be performed in probabilistic polynomial time \cite{VoightHilbert} by extracting square roots modulo $\frakp$.
\end{proof}

Proposition \ref{indefclassset} thus proves the indefinite case of Theorem B.  

\begin{rmk}
{\rm In practice, in Algorithm \ref{ClassSetEichler} one may wish to find representatives of $\Cl \calO$ with the smallest norm possible; one can then simply find small representatives $\fraka$ for each ideal class of $\Cl_{S_\infty} \Z_F$ (using the LLL algorithm \cite{LLL}, part of the algorithms used in the algorithm described in Proposition \ref{computeClZF}) and then repeat Step 1 for each such ideal $\fraka$.}
\end{rmk}

To solve Problem \ref{conjugacyclasses} (\textsf{ConjClassSet}) for Eichler orders, we amend Algorithm \ref{ClassSetEichler} as follows.  

\begin{alg} \label{ConjSetEichler}
Let $\calO \subset B$ be an indefinite Eichler order of discriminant $\frakd$.  This algorithm solves Problem \ref{conjugacyclasses} \textup{(\textsf{ConjClassSet})}.
{\rm \begin{enumerate}
\item Let $H$ be the subgroup of $\Cl_{S_\infty} \Z_F$ generated by $2 \Cl \Z_F$ and $[\frakp^e]$ for all $\frakp^e \parallel \frakd$ with $e$ odd.
\item For each $\fraka$ in a set of representatives for $\Cl_{S_\infty} \Z_F/H$, let $\alpha \in \calO$ be such that $\nrd(\alpha)$ is a uniformizer of $\fraka$ and let $I_\fraka := \alpha\calO+\fraka\calO$.
\item Return the set $\{\calO_L(I_\fraka)\}_{\fraka}$.
\end{enumerate}}
\end{alg}

\begin{proof}[Proof of correctness]
By Proposition \ref{alltwosidedclasses}, the image of the reduced norm of the set of two-sided ideal classes maps is exactly $H$.  It then follows from Proposition \ref{conjareequiv} that one recovers all conjugacy classes of Eichler orders of level $\frakN$ as $\calO_L(I_\fraka)$ for the right $\calO$-ideals $I_\fraka$ with $\fraka$ as in Step 2.
\end{proof}

Finally, we are left with Problem \ref{computeprincipal} (\textsf{IsPrincipal}).  Let $I$ be a right fractional $\calO$-ideal.  Again, by Eichler's theorem (Proposition \ref{ECcond}), we see that $I$ is principal if and only if $\nrd(I) \subset \Z_F$ is trivial in $\Cl_{S_\infty} \Z_F$, and the latter can be tested as in \S 2.  In other words, simply testing if a right $\calO$-ideal is principal is no harder than testing if an ideal is principal in $\Z_F$.  

To then actually exhibit a generator for a principal ideal, we rely upon the following standard lemma (see Pizer \cite[Proposition 1.18]{Pizer}).

\begin{lemma} \label{lemprin}
Let $I$ be a right invertible fractional $\calO$-ideal.  Then $\gamma \in I$ generates $I$ if and only if $\nrd(\gamma)\Z_F=\nrd(I)$.
\end{lemma}

% \begin{proof}
% Without loss of generality, assume that $I \subset \calO$ is integral.  If one first defines the norm $\N(I)$ of $I$ to be the product of the primes of $\calO$ occuring in a composition series for $\calO/I$ as a $\calO$-module (as in Reiner \cite[Equation 24.1]{Reiner}), then the statement $\gamma \calO = I$ if and only if $\N(\gamma \calO) = \N(I)$ is obvious.  Since $[B:F]=4$, we have $\N(I)=\nrd(I)^2$ \cite[Theorem 24.11]{Reiner}, and the result follows.
% \end{proof}

By Lemma \ref{lemprin}, the right ideal $I$ is principal if and only if there exists $\gamma \in I$ such $\nrd(\gamma)=\nrd(I)=c\Z_F$ (with $v(c) > 0$ for all $v \in S_\infty$).  Unfortunately, since $B$ is indefinite, the norm form $\Tr\nrd:B \to \Q$ is not positive definite, hence it does not induce the structure of a (definite) lattice on $I$ (in the definite case it will, see \S\S 6--7).  One option is to use a form of indefinite lattice reduction (as in Ivanyos-Sz\'ant\'o \cite{IS}).  We instead find a substitute quadratic form which will still allow us to find ``small'' elements.  When $F$ is totally real and $B$ has a unique split real place, such a form has been found \cite{Voight-funddom}, and inspired by this result we make the following definitions.

Let $B=\quat{a,b}{F}$.  For an infinite place $v$ of $F$ and $\gamma=x+y\alpha+z\beta+w\alpha\beta$, define
\begin{equation} \label{defred}
Q_v(\gamma) = |v(x)|^2+|v(a)||v(y)|^2 + |v(b)||v(z)|^2 + |v(ab)||v(w)|^2.
\end{equation}
We then define the \emph{absolute reduced norm} by
\begin{align*}
Q:B &\to \R \\
\gamma &\mapsto \textstyle{\sum_v} Q_v(\gamma);
\end{align*} 
by construction, the form $Q$ is positive definite and gives $I$ the structure of a definite $\Z$-lattice of rank $4[F:\Q]$.

\begin{rmk}
{\rm The form $Q$ is clearly only one of many choices for such a positive definite form, and so one may reasonably try to understand what the cone of such forms corresponds to.

When $F$ is totally real and $B$ has a unique split real place, the choice of the positive definite quadratic form corresponds to the choice of a center $p$ for a Dirichlet fundamental domain in the upper half-plane $\frakH$ and at the split place measures the \emph{inverse radius} of the corresponding isometric circle \cite{Voight-funddom}.  The same is true for a quaternion algebra of arbitrary signature as follows.  If $B$ has $g$ split real places and $s$ (split) complex places, then the group $\calO_1^*$ of units of $\calO$ of reduced norm $1$ embeds in $SL_2(\R)^g \times SL_2(\C)^s$ and acts on $\frakH^g \times (H^3)^s$ discretely, where $\frakH$ (resp.\ $H^3$) denotes the upper half-plane (resp.\ hyperbolic $3$-space) (see e.g.\ Beardon \cite{Beardon} and Elstrodt, et al.\ \cite{Elstrodt}).  In this case, the choice of a positive definite quadratic form corresponds again to the choice of a center $p$ for a Dirichlet domain and at each place measures an inverse radius, either of an isometric circle or sphere.  The above choice of form $Q$ corresponds to a (normalized) choice of center $p=(i,\dots,i,j,\dots,j)$.  

Because of the connection with the classical theory of positive definite quadratic forms on a real quadratic field, which can be understood more generally from the perspective of Arakelov theory \cite{Schoof}, we view these observations as the beginning of a form of noncommutative Arakelov theory and leave it as a subject for further investigation.}
\end{rmk}

We then have the following algorithm.

\begin{alg} \label{findprinc}
Let $I \subset \calO$ be a right fractional $\calO$-ideal.  This algorithm solves Problem \textup{\ref{computeprincipal}} \textup{(\textsf{IsPrincipal})}.
{\rm 
\begin{enumerate}
\item Compute $\nrd(I) \subset \Z_F$ and test if $\nrd(I)$ is principal; if not, then return \textup{\textsf{false}}.  Otherwise, let $\nrd(I)=c\Z_F$.
\item Determine if there exists a unit $u \in \Z_F^*$ such that $v(uc) > 0$ for all ramified (real) places $v$; if not, then return \textup{\textsf{false}}.  Otherwise, let $c := uc$ and initialize $\alpha := 1$.
\item If $c \in \Z_F^*$ , return $\alpha$.  Otherwise, view $I$ as a $\Z$-lattice equipped with the quadratic form $Q$.  Reduce $I$ using the LLL algorithm \cite{LLL}.  By enumerating short elements in $I$, find $\gamma \in I$ such that $\nrd(\gamma)=cd$ with $\N(d)<\N(c)$.  Let $\alpha := \gamma\alpha/d$, let $I := d\gamma^{-1} I$, and let $c := d$, and return to Step $2$.
\end{enumerate}}
\end{alg}

\begin{proof}
In Step 2, we have $\nrd(d\gamma^{-1}I)=d^2/(cd)\nrd(I)=d\Z_F$, and so the algorithm terminates since in each step $\N\nrd(I) \in \Z_{>0}$ decreases.  The algorithm gives correct output since $d\gamma^{-1} I=\alpha \calO$ if and only if $I=(\gamma\alpha/d)\calO$.
\end{proof}

% \begin{rmk}
% One may also improve this algorithm by generalizing the algorithm of Demb\'el\'e-Donnelly \cite[\S 2.2]{DembeleDonnelly} to rescale the lattice at the real ramified places (see Algorithm \ref{findprincdef}).
% \end{rmk}

In practice, Algorithm \ref{findprinc} runs quite efficiently and substantially improves upon a more na\"ive enumeration.  However, we are unable to prove any rigorous time bounds for Algorithm \ref{findprinc}.  Already the first step of the algorithm requires the computation of the class group $\Cl \Z_F$; and even if we suppose that the class group has been precomputed, there do not appear to be rigorous time bounds for the principal ideal testing algorithm \cite[Algorithm 6.5.10]{Cohen} (see \S 1).  With that proviso, given the generator $c$ as in Step 1, we can measure the value of the LLL-step as follows.

\begin{lemma}
There exists $C(\calO) \in \R_{>0}$, depending on $\calO$, such that for every principal fractional ideal $I$ of $\calO$, the first basis element $\gamma$ in the $LLL$-reduced basis of Algorithm \textup{\ref{findprinc}} satisfies
\[ |\N(\nrd(\gamma))| \leq C(\calO) \N(\nrd(I)). \]
\end{lemma}

% Choose a representation $B=\quat{a,b}{F}$ such that the order $\Lambda$ generated by the standard generators $\alpha,\beta$ has $\Lambda \subset \calO$ as in (\ref{quateq}).

%Let $L$ denote the lattice $I$ with the norm $N$.  We first claim that
%\[ \det(L)=\frac{d_F^2 \N(\nrd(I))^2 |\N(ab)|^2}{[\calO:\Lambda]}. \]  

\begin{proof}
Suppose that $I=\xi \calO$.  The $F$-endomorphism of $B$ given by left multiplication by $\xi$ has determinant $\nrd(\xi)^2$, and it follows that the corresponding $\Q$-endomorphism of $L \otimes_\Z \Q$ has determinant $\N(\nrd(\xi))^2$.  Hence
\[ \det(I)=\det(\xi \calO)=\N(\nrd(\xi))^4 \det(\calO). \]
% \[ \det(I)=\det(\xi \calO)=\N(\nrd(\xi))^4 \det(\calO)= \N(\nrd(I))^4 \frac{\det(\Lambda)}{[\calO:\Lambda]^2} \]
% since $\Lambda \subset \calO$. % Now it is easy to see that $\det(\Lambda)=d_F^2 |\N(ab)|^2$, and the claim follows.

Now, for any $\gamma=x+y\alpha+z\beta+w\alpha\beta$, from (\ref{defred}) we have
\[ |v(\nrd(\gamma))| \leq Q_v(\gamma) \leq Q(\gamma) \]
for all places $v$.  Thus, the output of the LLL algorithm \cite[Proposition 1.9]{LLL} yields $\gamma \in I$ which satisfies
\[ |v(\nrd(\gamma))| \leq Q(\gamma) \leq 2^{(4n-1)/4} \det(I)^{1/(4n)}=2^{(4n-1)/4} \det(\calO)^{1/(4n)}
\N(\nrd(I))^{1/n}. \]
We conclude that
\[ |\N(\nrd(\gamma))| = \prod_v |v(\nrd(\gamma))| \leq 2^{(4n^2-n)/4}\det(\calO)^{1/4}\N(\nrd(I)) \]
as claimed.
\end{proof}

From Lemma \ref{lemprin}, we conclude that the algorithm produces elements which are close to being generators.

\section{Computing the class number for definite orders}

In this section, we discuss the difficulty of solving Problem \ref{computeclassnumber} (\textsf{ClassNumber}) in the definite case.  Throughout this section, let $B$ denote a totally definite quaternion algebra of discriminant $\frakD$.  Here, the class number is governed by the Eichler mass formula.  

Given an ideal $\frakN$ (coprime to $\frakD$), the \emph{mass} is defined to be the function
\begin{equation} \label{massdef}
M(\frakD,\frakN)=2^{1-n}|\zeta_F(-1)|h(\Z_F)\Phi(\frakD)\Psi(\frakN)
\end{equation}
where
\begin{equation}
\Phi(\frakD) = \prod_{\frakp \mid \frakD}\left(\N(\frakp)-1\right) \quad \text{and} \quad
\Psi(\frakN) = \N(\frakN)\prod_{\frakp \mid \frakN}\left(1+\frac{1}{\N(\frakp)}\right).
\end{equation}
The mass of an Eichler order $\calO \subset B$ of level $\frakN$ is defined to be $M(\calO)=M(\frakD,\frakN)$.

The class number of an Eichler order differs from its mass by a correction factor coming from torsion, as follows.  An \emph{embedded elliptic subgroup} in $B$ is an embedding $\mu_q \hookrightarrow \calO^*/\Z_F^*$, where $q \in \Z_{\geq 2}$ and $\calO$ is an Eichler order, such that the image is a maximal (cyclic) subgroup of $\calO^*/\Z_F^*$; the \emph{level} of the embedding is the level of $\calO$.  An \emph{elliptic cycle} is a $B^*$-conjugacy class of embedded elliptic subgroups.  Let $e_q(\frakD,\frakN)$ denote the number of elliptic cycles of $B^*$ of order $q$ and level $\frakN$.

\begin{proposition}[Eichler mass formula {\cite[Corollaire V.2.5]{Vigneras}}] \label{massformula}
Let $\calO \subset B$ be an Eichler order of level $\frakN$.  Then
\[ h(\calO) = M(\frakD,\frakN) + \sum_q e_q(\frakD,\frakN)\left(1-\frac{1}{q}\right). \]
\end{proposition}

\begin{rmk} \label{altmassformula}
{\rm A variant of the Eichler mass formula \cite[Corollaire V.2.3]{Vigneras} which is also useful for algorithmic purposes (see Remark \ref{usealtmassformula}) reads
\[ M(\frakD,\frakN) = \sum_{[I] \in \Cl \calO} \frac{1}{[\calO_L(I)^*:\Z_F^*]}. \]}
\end{rmk}

We first characterize the embedding numbers $e_q(\frakD,\frakN)$.  Given an embedded elliptic subgroup $\mu_q \hookrightarrow \calO^*/\Z_F^*$ of level $\frakN$, the image of $\mu_q$ generates a quadratic subring $R \subset \calO$; such an embedding $R \hookrightarrow \calO$ with $RF \cap \calO=R$ is said to be an \emph{optimal embedding}.  Conversely, to every optimal embedding $\iota:R \hookrightarrow \calO$, where $R$ is a quadratic $\Z_F$-order with $[R^*:\Z_F^*]=q$ and $\calO$ is an Eichler order of level $\frakN$, we have the embedded elliptic subgroup $R_{\text{tors}}^*/\Z_F^* \cong \mu_q \hookrightarrow \calO^*$.  This yields a bijection
\begin{center}
$\{\text{Elliptic cycles of $B^*$ of order $q$ and level $\frakN$}\}$ \\
$\updownarrow$ \\
$\biggl\{$
\parbox{0.66\textwidth}{$B^*$-conjugacy classes of optimal embeddings $\iota:R \hookrightarrow \calO$ \\ \phantom{xx}with $[R^*:\Z_F^*]=q$ and $\calO$ an Eichler order of level $\frakN$} $\biggr\}$.
\end{center}

The quadratic $\Z_F$-orders $R$ with $[R^*:\Z_F^*]=q$ come in two types.  Either we have $R_{\text{tors}}^* \cong \mu_{2q}$ and we say $R$ is \emph{fully elliptic}, or $[R^*:\Z_F^*R_{\text{tors}}^*]=2$ and we say $R$ is \emph{half elliptic}.  We see that if $R$ is half elliptic then in particular $R \subset \Z_F[\sqrt{-\eps}]$ for $\eps$ a totally positive unit of $\Z_F$.

The (global) embedding numbers $e_q(\frakD,\frakN)$ can then be computed by comparison to the local embedding numbers
\[ m(R_\frakp,\calO_\frakp)=\#\{\text{$\calO_\frakp^*$-conjugacy classes of optimal embeddings $\iota:R_\frakp \hookrightarrow \calO_\frakp$}\} \]
where $\calO_\frakp$ is a $\frakp$-local Eichler order of level $\frakN$.

\begin{lemma}[{\cite[p.~143]{Vigneras}}] \label{embedclassno}
We have
\[ e_q(\frakD,\frakN)=\frac{1}{2} \sum_{[R^*:\Z_F^*]=q} h(R) \prod_\frakp m(R_\frakp,\calO_\frakp). \]
\end{lemma}

There are formulas \cite[\S 2]{VoightFuchs} for the number of local embeddings $m(R_\frakp,\calO_\frakp)$, for example:
\begin{equation} \label{woot}
m(R_\frakp,\calO_\frakp)=
\begin{cases}
1, & \text{if $\frakp \nmid \frakD\frakN$;} \\
\displaystyle{1-\legen{K_q}{\frakp}}, & \text{if $\frakp \mid \frakD$ and $\frakp \nmid \frakf(R)$;} \\
\displaystyle{1+\legen{K_q}{\frakp}}, & \text{if $\frakp \parallel \frakN$;}
\end{cases}
\end{equation}
here, we let $K_q=F(\zeta_{2q})$ and $\frakf(R)$ denotes the conductor of $R$ (in $\Z_{K_q}$).  In particular, we have by Equation (\ref{woot}) that $m(R_\frakp,\calO_\frakp)=1$ for almost all $\frakp$.  

We now discuss the computability of the terms in the formula of Proposition \ref{massformula}.  To compute the mass, we will use the following proposition.

\begin{proposition} \label{computezeta}
The value $\zeta_F(-1) \in \Q$ can be computed using $O\bigl(d_F^{3/2}\log^4 d_F\bigr)$ bit operations.
\end{proposition}

\begin{proof}
From the functional equation for the Dedekind zeta function, we have
\begin{equation} \label{zetaFm1}
\zeta_F(-1)=\left(\frac{-1}{2\pi^2}\right)^n d_F^{3/2} \zeta_F(2).
\end{equation}
From (\ref{massdef}) and Proposition \ref{massformula}, we have $\zeta_F(-1) \in \Q$, in fact, $\zeta_F(-1)$ has denominator bounded by $Q$, the least common multiple of all $q \in \Z_{\geq 2}$ such that $[F(\zeta_{2q}):F] = 2$.  

We compute an approximation to $\zeta_F(-1)$ from the Euler product expansion for $\zeta_F(2)$, as follows (see also Buchmann-Williams \cite[\S 2]{BuchWill}, or Dokchitser \cite{Dokchitser} for a more general approach).  For $P \in \Z_{\geq 2}$, let
\[ \zeta_{F,\leq P}(s)=\prod_{\N\frakp \leq P} \left(1-\frac{1}{\N\frakp^s}\right)^{-1} \]
denote the truncated Euler product for $\zeta_F(s)$, where we take the product over all primes $\frakp$ of $\Z_F$ for which $\N\frakp \leq P$.  Note that for $s > 1$ real we have
\[ \zeta_{F,\leq P}(s) \leq \zeta_F(s) = \prod_{\frakp} \left(1-\frac{1}{\N\frakp^s}\right)^{-1} \leq \prod_p \left(1-\frac{1}{p^s}\right)^{-n} = \zeta(s)^n. \]
Now we estimate
\begin{align*}
0 < \frac{\zeta_F(2)}{\zeta_{F,\leq P}(2)}- 1 &= \prod_{\N\frakp > P} \left(1-\frac{1}{\N\frakp^2}\right)^{-1} - 1 \\
&= \sum_{\N\frakp >P} \frac{1}{N\frakp^2} + 
\sum_{\N\frakp \geq \N\frakq>P} \frac{1}{N(\frakp\frakq)^2} + \dots  \\
& \leq \sum_{p>P} \frac{n}{p^2} + \sum_{p,q>P} \frac{n^2}{(pq)^2} + \dots 
\leq \sum_{x>P} \frac{n}{x^2} + \sum_{x>P^2} \frac{n^2}{x^2} + \dots \\
& \leq \frac{n}{P} + \frac{n^2}{P^2} + \dots = \frac{1}{(P/n)-1}.
\end{align*}
It follows that $\zeta_F(2) - \zeta_{F,\leq P}(2) < \eps$ whenever $P>n(1+\zeta_{F,\leq P}(2)/\eps)$ which is satisfied when
\[ P>n\left(1+\frac{\zeta(2)^n}{\eps}\right) = n\left(1+\frac{(\pi^2/6)^n}{\eps}\right).  \]
To obtain the value $\zeta_F(-1)$ within an interval of length smaller than $1/(2Q)$, from (\ref{zetaFm1}) it suffices to take
\[ \eps=\frac{1}{2Q} \frac{(2\pi^2)^n}{d_F^{3/2}}. \]

To estimate $Q$, we note that if $[F(\zeta_{2q}):F]=2$ then in particular $F$ contains the totally real subfield $\Q(\zeta_{2q})^+$ of $\Q(\zeta_{2q})$, hence $\phi(q)/2 \mid n$.  Since the fields $\Q(\zeta_{2q})^+$ are linearly disjoint for $q$ a power of a prime, we have $Q=O(n)$.  (We note this bound is best possible in terms of $n$, since after all we may take $F=\Q(\zeta_{2q})^+$.)  Putting these together, we need to evaluate the truncated Euler product with 
\[ P=O\left(n \left(\frac{\pi^2}{6}\right)^n\frac{1}{\eps}\right)=O\left(\frac{n^2}{12^n} d_F^{3/2}\right). \]

Evaluating $\zeta_{F,\leq P}(2)$ amounts to factoring a degree $n$ polynomial over $\F_p$ for all primes $p \leq P$; each such factorization can be performed using a repeated squaring operation, requiring $O(n^3 \log p)$ operations in $\F_p$ (see e.g.\ the survey by von zur Gathen and Panario \cite{Gathen}) so time $O(n^3 \log^3 p)$, hence altogether time $O(n^3 P \log^2 P)$ by the prime number theorem, so the computation requires
\[ O\left(\frac{n^5}{12^n}d_F^{3/2}\log^3 d_F\right)=O(d_F^{3/2}\log^3 d_F) \]
operations with real numbers of precision $O(1/\eps)=O(d_F^{3/2})$, requiring therefore $O(d_F^{3/2} \log^4 d_F)$ bit operations.
\end{proof}

\begin{rmk}
{\rm If $F$ is an abelian field, then $\zeta_F(-1)$ can be computed much more efficiently in terms of Bernoulli numbers \cite{Washington}.}
\end{rmk}

Putting these pieces together, we now prove the following theorem.

\begin{theorem} \label{existstotdef}
There exists a probabilistic algorithm which, given an Eichler order $\calO$ in a definite quaternion algebra $B$ with factored discriminant $\frakd$, solves Problem \textup{(\textsf{ClassNumber})} in time
\[ O\bigl(d_F^{3/2}\log^4 d_F+\log^2 \N \frakd \bigr) \] 
and the time to solve $O(2^n)$ instances of Problem \textup{\textsf{(ClassUnitGroup)}} with fields of discriminant of size $O(d_F^{5/2})$.
\end{theorem}

\begin{proof}
First, we compute the factored discriminant $\frakD$ of $B$ and level $\frakN$ of $\calO$ by computing Hilbert symbols \cite{VoightHilbert}: given the factorization of the discriminant $\frakd$ of $\calO$, for each $\frakp \mid \frakd$ one can determine whether $\frakp \mid \frakD$ or $\frakp \mid \frakN$ in deterministic time $O(\log^2 \N \frakp)$.

We compute $h(\calO)$ from the Eichler mass formula (Proposition \ref{massformula}), with the mass $M(\frakD,\frakN)$ given as in Equation (\ref{massdef}).  Given the factorization of $\frakD$ and $\frakN$, we can compute $\Phi(\frakD)$ and $\Psi(\frakN)$ in time $O(\log^2 \N(\frakD\frakN))$.  One recovers $h(\Z_F)$ from the given algorithm to solve Problem \textup{\textsf{(ClassUnitGroup)}}.  By Proposition \ref{computezeta}, we can compute $|\zeta_F(-1)|$ in time $O(d_F^{3/2}\log^4 d_F)$.

We now estimate the time to compute the correction term $\sum_q e_q(\frakD,\frakN)(1-1/q)$ in Lemma \ref{embedclassno}, up to further instances of Problem \textup{\textsf{(ClassUnitGroup)}}. 

As in the proof of Proposition \ref{computezeta}, if $e_q(\frakD,\frakN) \neq 0$ then $q=O(n)$.  Thus, by factoring $n$ (which can be done in negligible time using trial division) we obtain a finite set of $O(n)$ possible values of $q$.  For each such $q$, we can factor the $q$th cyclotomic polynomial over $F$ in deterministic time $(qn\log d_F)^{O(1)}$ (see Lenstra \cite{AKLenstra}) and determine if each of its irreducible factors has degree $2$.  (In practice, one might further restrict the set of possible $q$ by verifying that $q \mid d_F$ if $q \geq 5$ and that for sufficiently many small prime ideals $\frakp$ of $\Z_F$ coprime to $qd_F$ that $q \mid (N\frakp^2-1)$, but this makes no difference in the analysis of the runtime.)  Since $q=O(n)$ and $n=O(\log d_F)$ by the Odlyzko bounds \cite{Odlyzko}, this becomes simply $(\log d_F)^{O(1)}$.

We pause to compute an estimate of discriminants.  Let $q \in \Z_{\geq 2}$ be such that $K=F(\zeta_{2q})$ is quadratic over $F$.  We estimate the discriminant $d_{K}=d_F^2 \N_{F/\Q}\frakd_{K/F}$.  Let $\lambda_{2q}=\zeta_{2q}+1/\zeta_{2q}$; then we have $\Q(\zeta_{2q})^+=\Q(\lambda_{2q})$ and $\zeta_{2q}^2-\lambda_{2q}\zeta_{2q}+1=0$, so $\frakd_{K/F} \mid (\lambda_{2q}^2-4)\Z_F$.  Therefore
\begin{align*} 
\N_{\Q(\zeta_{2q})^+/\Q}(\lambda_{2q}^2-4) &= \N_{\Q(\zeta_{2q})^+/\Q}(\zeta_{2q}-1/\zeta_{2q})^2 = \N_{\Q(\zeta_{2q})/\Q}(\zeta_q-1) \\
&= \begin{cases}
p, &\text{ if $q=p^r$ is an odd prime power}; \\
4, &\text{ if $q=2^r$}; \\
1, &\text{ otherwise}.
\end{cases}
\end{align*}
So $d_{K}=d_F^2 \N_{F/\Q} \frakd_{K/F} \leq d_F^2 \N_{F/\Q} (\lambda_{2q}^2-4) = d_F^2 p^{2n/\phi(2q)}$ if $q=p^r$ is an odd prime power, and similarly $d_{K} \leq d_F^2 4^{2n/\phi(2q)}$ if $q=2^r$ and $d_{K}=d_F^2$ if $q$ is not a power of a prime.  For $q$ a power of a prime $p$, by the conductor-discriminant formula \cite[Theorem 3.11]{Washington} we have
$p^{\phi(2q)/2-1} \mid d_{\Q(\zeta_{2q})^+}$, so since $d_{\Q(\zeta_{2q})^+}^{2n/\phi(2q)} \mid d_F$ we have $p^n = O(d_F)$.  Thus we have $p^{2n/\phi(2q)}=O(d_F^{2/\phi(2q)})=O(d_F^{1/2})$ when $q=5$ or $q \geq 7$, and hence $d_{K}=O(d_F^{5/2})$; but for $q \leq 4$ this also holds, since by the Odlyzko bounds there are only finitely many (totally real) number fields with $d_F \leq 16^n$.  When $K=F(\sqrt{-\eps})$ for $\eps$ a totally positive unit of $\Z_F$, we have $d_K = O(d_F^2 4^n) = O(d_F^{5/2})$ as well.  Thus in all cases, $d_K=O(d_F^{5/2})$.   

To compute $e_q(\frakD,\frakN)$, we need to compute $h(R)$ for all $R \subset K$ with $[R^*:\Z_F^*]=q$.  First, suppose that $R$ is fully elliptic.  Then $\Z_F[\zeta_{2q}] \subset R \subset \Z_{K}$.  We have the bound $\frakf(R)^2 \mid (\lambda_{2q}^2-4) \mid p$ on the conductor of $R$ as above, which implies there are at most $O(2^{n/\phi(2q)})=O(d_F)$ such orders, corresponding to the possible products of ramified primes over $p$ in the extension $F/\Q(\zeta_{2q})^+$.  Each can be constructed by computing $\frakp$-overorders for $\frakp^2 \mid (\lambda_{2q}^2-4) \mid p$, requiring time $q^{O(1)}$.  Now if $\frakf=\frakf(R)$, we have
\[ h(R) = \frac{h(\Z_{K_q})}{[\Z_{K_q}^*:R^*]} \frac{\#(\Z_K/\frakf\Z_K)^*}{\#(\Z_F/\frakf\Z_F)^*} \]
and the time required to compute these terms is negligible except for $h(\Z_{K_q})$, for which we may call our hypothesized algorithm.  

% Although there are $O(n)$ values of $q$, there are $O(\log n)$ such $q$ which are prime powers, and we claim that it is enough to compute $\Z_{K_q}$ for $q$ a prime power: indeed, if $[F(\zeta_{2q}):F]=2$, then immediately $[F(\zeta_{2d}):F]=2$ for all $d \mid q$ with $d \in \Z_{\geq 2}$ so $F(\zeta_{2d})=F(\zeta_{2q})$ for such $d$.

Next, suppose that $R$ is half elliptic.  In this case, we have $\Z_F[\sqrt{-\eps}] \subset R \subset \Z_K$.  Here, we have the bound $\frakf(R)^2 \mid 4$, hence there are at most $O(2^n)=O(d_F)$ such orders, and the arguments in the preceding paragraph apply.

In all, we have at most $O(2^n)$ imaginary quadratic fields $K$ to consider, the worst case being when $F$ has a fundamental system of units which are totally positive. 

Finally, the calculation of the local embedding number $m(R_\frakp,\calO_\frakp)$ can be accomplished in time $O(\log^2 (\N\frakp))$ for $\frakp$ odd \cite[Proposition 2.5]{VoightFuchs} by computing a Legendre symbol (and in time $O(1)$ after that) and in time $n^{O(1)}=(\log d_F)^{O(1)}$ for $\frakp$ even \cite[Remark 2.6]{VoightFuchs}.

In all, aside from the time to compute class numbers, since there are $O(n)$ values of $q$, we can compute the correction term in time 
\[ O(\log^2 \N(\frakD\frakN)) + d_F(\log d_F)^{O(1)} + (\log d_F)^{O(1)}. \]  
Combining this estimate with the time to compute the mass, the result follows.
\end{proof}

\begin{corollary} \label{defprob}
For a fixed number field $F$, Problem \textup{\ref{computeclassnumber} (\textsf{ClassNumber})} for definite orders $\calO$ is probabilistic polynomial-time reducible to the factorization of the discriminant of $\calO$.
\end{corollary}

\begin{proof}
Only the quantities $\Phi(\frakD)$, $\Psi(\frakN)$, and $m(R_\frakp,\calO_\frakp)$ depend on $\calO$---the others can be precomputed for fixed $F$.  These quantities can be computed in probabilistic polynomial time knowing the factorization of $\frakD$ and $\frakN$.
\end{proof}

Putting together Corollary \ref{indefO(1)}, Theorem \ref{existstotdef}, and Corollary \ref{defprob}, we have proven the main Theorem A and its corollary.  

We conclude this section by discussing the role of factoring the ideals $\frakD, \frakN$.  It is well known that factoring ideals in the ring $\Z_F$ is probabilistic polynomial-time equivalent to factoring integers.  But already for imaginary quadratic fields, it is also well known that an algorithm for (\textsf{ClassUnitGroup}) can be employed to factor integers (two such ``class group'' methods are attributed to Schnorr-Lenstra-Shanks-Pollard-Atkin-Rickert and Schnorr-Seysen-Lenstra-Lenstra-Pomerance).  A noncommutative analogue of this result, already in the simplest case where $F=\Q$, is the following.

\begin{proposition} \label{factorints}
The problem of factoring integers $a$ with $\omega(a)=O(\log \log a)$ prime factors is probabilistic polynomial-time reducible to Problem \textup{\textsf{(ClassNumber)}} for definite quaternion algebras over $\Q$.
\end{proposition}

\begin{proof}
For $F=\Q$, the class number is given simply by 
\[ h(\calO)=\frac{1}{12}\phi(D)\psi(N) + \frac{e_2(D,N)}{2} + \frac{e_3(D,N)}{3}. \]
Let $a \in \Z_{>0}$ be an integer to be factored, which we may assume has $\gcd(a,26)=1$ and is not a prime power.

We consider quaternion algebras of the form $B=\quat{-ac,-13b}{\Q}$ where $b,c \in \Z_{>0}$ are chosen as follows.  We choose $c \in \{1,2\}$ so that $(-ac/13)=-1$: therefore the algebra $B$ is ramified at $13$ so $13 \mid D$, and it follows that $e_2(D,N)=e_3(D,N)=0$ by (\ref{woot}), hence
\[ h(\calO)=\phi(D/13)\psi(N). \]
For simplicity, we assume $c=1$; the same argument applies when $c=2$.  

We choose $b$ to be a random squarefree factored positive integer modulo $4a$ (see e.g.\ Bach \cite{Bach}) with $\gcd(b,26)=1$.  If $\gcd(a,b) \neq 1$ and $a \neq 1,b$, we have factored $a$.  Otherwise, with probability at least $(1/2)^{\omega(a)}$, we have $(-13b/p)=1$ for all $p \mid a$, so then $p \nmid D$ for all $p \mid a$; since $\omega(a)=O(\log\log a)$ by assumption, after $O(\log a)$ attempts this condition will be satisfied with probability at least $1/2$.

Now compute an order $\calO \subset B$ which is locally generated by the standard generators $\alpha,\beta$ for all primes $p \neq 2$ and which is maximal at $2$.  Then $\calO$ is an order with discriminant $13ab\eps=DN$ where $\eps=1$ or $2$.  We claim that $\calO$ is Eichler.  Since $\gcd(a,b)=1$ and $b$ is squarefree, and further $a,b$ are odd, it suffices to show that $\calO$ is Eichler at $p \mid a$.  But for such $p^e \parallel a$, by assumption we have $(-13b/p)=1$ so $B$ is unramified at $p$, and there exists $t \in \Z_p^*$ such that $t^2=-13b$.  Then the embedding $\calO \hookrightarrow \calO \otimes_{\Z} \Z_p$ by
\[ \alpha \mapsto \begin{pmatrix} 0 & 1 \\ ac & 0 \end{pmatrix}, \quad \beta \mapsto \begin{pmatrix} t & 0 \\ 0 & -t \end{pmatrix} \]
realizes $\calO$ as an Eichler order of level $p^e$.  

It follows that
\[ h(\calO)=\prod_{p \mid (D/13)}(p-1) \prod_{p^e \parallel N} p^{e-1}(p+1); \]
if $a$ is not squarefree then $\gcd(h(\calO),a)$ yields a prime factor of $a$, so we may assume that $a$ is squarefree.  

We now show how to recover a prime factor of $a$ given the values $h(\calO)$.  For each prime $p \mid b$, by computing $(-ac/p)$ we can determine if $p$ is ramified in $B$ or not, accordingly contributing a factor $p-1$ or $p+1$ to $h(\calO)$.  Dividing $h(\calO)$ by these factors, we may compute the integer
\[ g(b)=\prod_{p \mid D'}(p-1) \prod_{p \mid N'} (p+1) \]
and $\eps a=D'N'$.  We compute one such value $g(b)$ and then many more values $g(b')$, and we claim that $\gcd(g(b)+g(b'),a)$ will find a factor of $a$.  Indeed, with probability at least $\omega(a)(1/2)^{\omega(a)+1}$, we will have the following:
\begin{center}
There exists a prime $q \mid a$ for which $(-13b/q)=-(-13b'/q)$ and $(-13b/p)=(-13b'/p)$ for all $p \mid \eps a$ with $p \neq q$.
\end{center}
It follows then that $g(b) \equiv -g(b') \pmod{q}$ as claimed.  Again, since by assumption $\omega(a)=O(\log \log a)$, after $O(\log a)$ attempts we will factor $a$ with probability at least $1/2$.
\end{proof}

\section{Ideal principalization for definite orders}

In this section, we discuss the totally definite case of Problem \ref{computeprincipal} (\textsf{IsPrincipal}), which by Lemma \ref{reductionisomtest} is equivalent to Problem \ref{computeisomorphism} (\textsf{IsIsomorphic}).  We exhibit an algorithm (Algorithm \ref{findprincdef}) to solve this problem and analyze its running time (Proposition \ref{dembdonn}). This algorithm will be used in the next section to solve Problem \ref{computeidealclasses} (\textsf{ClassSet}) for definite orders.  Throughout this section, let $B$ be a totally definite quaternion algebra over a totally real field $F$ of discriminant $\frakD$ and let $\calO \subset B$ be an Eichler order of level $\frakN$.  

Our algorithm is similar to the indefinite case (Algorithm \ref{findprinc}), where now we are in the easier situation that $\Tr\nrd:B \to \Q$ is a positive definite quadratic form on the $\Q$-vector space $B$.  We prove below that we can reduce the principalization problem to a shortest lattice vector problem (based initially on an idea due to Demb\'el\'e and Donnelly \cite[\S 2.2]{DembeleDonnelly}.)  

First of all, if $I$ is a principal right fractional $\calO$-ideal, then $\nrd I=c\Z_F$ is a principal ideal of $\Z_F$.  To compute such a generator $c \in \Z_F$ over a general number field $F$, we refer to the discussion following Proposition \ref{computeClZF}.  For $F$ totally real as in this section, we first show that a principal ideal has a generator of polynomial size, with the constant depending on $F$.  For a real place $v$ of $F$, corresponding to an embedding $v:F \hookrightarrow \R$, we abbreviate $a_v=v(a)$ for $a \in F$.

\begin{lemma} \label{princZF}
For a principal fractional ideal $\fraka \subset F$, there exists a generator $a$ for $\fraka$ with $|a_v|=O((\N \fraka)^{1/n})$ for all real places $v$ of $F$, where the implied constant is effectively computable and depends on the field $F$.  
\end{lemma}

\begin{proof}
Consider the (Minkowksi unit) map 
\begin{align*}
\sigma:F &\to \R^n \\ 
u &\mapsto (\log |u_v|- (1/n)\log \N(u))_v
\end{align*} 
where $v$ runs over the real places of $F$.  The image $\sigma(F)$ is contained in the hyperplane $H:\sum_v x_v = 0$, and the image of the units $\sigma(\Z_F^*)$ forms a lattice (of full rank $n-1$) in $H$.  But since $\fraka$ is principal, the image $\sigma(\fraka) \subset H$ is simply a translate of this lattice.  It follows that there is a generator $a$ for $\fraka$ with $\sigma(a)$ inside a fundamental parallelopiped for $\sigma(\Z_F^*)$.  Therefore, if $u_1,\dots,u_{n-1} \in \Z_F^*$ are a basis for $\Z_F^*/\{\pm 1\}$, then for each real place $v$ we have 
\[ \left|\log |a_v|-\frac{1}{n}\log \N(a) \right|\leq \frac{1}{2}\sum_{i=1}^{n-1} |\log |u_{i,v}||. \]
Let $U_v = \prod_{i=1}^{n-1} e^{|\log |u_{i,v}||/2}$.  Then
\[ \frac{1}{U_v} \N(\fraka)^{1/n} \leq |a_v| \leq U_v \N(\fraka)^{1/n} \]
so if $U = \max U_v$ then
\[ \frac{1}{U_F} \N(\fraka)^{1/n} \leq |a_v| \leq U \N(\fraka)^{1/n} \]
for all $v$, as claimed.
\end{proof}

\begin{rmk} \label{princZFrmk}
{\rm Although we will not make use of it, we remark that Lemma \ref{princZF} gives a method for testing if an ideal $\fraka$ is principal over $\Z_F$ and, if so, exhibiting a generator.  We simply enumerate all elements $a \in \Z_F$ with $|a_v| \leq U \N(\fraka)^{1/n}$ for all $v$, where $U$ is as in the above proof; we can accomplish this by embedding $\Z_F \hookrightarrow \R^n$ by $a \mapsto (a_v)_v$ as a lattice (by the Minkowski embedding) and enumerating all lattice elements of (square) norm bounded by $nU^2 \N(\fraka)^{2/n}$.  By the analysis of Fincke and Pohst \cite[(2.12), (3.1)]{FinckePohst} (or alternatively Elkies \cite[Lemma 1]{ElkiesANTSIV}), we can enumerate all such elements in time $O( \N(\fraka) )$, where the implied constant depends only on $F$.}
\end{rmk}

We now present the main algorithm in this section.

\begin{alg} \label{findprincdef}
Let $I \subset \calO$ be a right fractional $\calO$-ideal and let $c \in \Z_F$ such that $\nrd I=c\Z_F$.  This algorithm solves Problem $2.4$ \textup{(\textsf{IsPrincipal})}.
{\rm \begin{enumerate}
\item Determine if there exists a unit $u \in \Z_F^*$ such that $c_v u_v > 0$ for all real places $v$; if not, then return \textup{\textsf{false}}.
\item For each totally positive unit $z \in \Z_{F,+}^*/\Z_F^{*2}$:
\begin{enumalgalph}
\item Let $\xi$ be a shortest vector of the $\Z$-lattice $I$ with respect to the rational quadratic form $\varphi_{(ucz)^{-1}}$.
\item If $\varphi_{(ucz)^{-1}}(\xi) = [F : \Q]$ then return \textup{\textsf{true}} and the element $\xi$.
\end{enumalgalph}
\item Return \textup{\textsf{false}}.
\end{enumerate}}
\end{alg}

\begin{rmk}
{\rm Note that if $F=\Q$ then in Step $2$ we have $z=u=1$. Hence the algorithm simply amounts to find a shortest vector in the $\Z$-lattice $I$ (with respect to the reduced norm form).}
\end{rmk}

\begin{proof}[Proof of correctness]
If $I$ is principal, then $\nrd I$ is generated by a totally positive element $uc$ where $u \in \Z_F^*$.  Then Lemma 4.8 implies that $\xi \in I$ generates $I$ if and only if $\nrd \xi = ucz$ for some $z \in \Z_{F,+}^*$.  To find such an element $\xi$, we only need to search for elements of norm $ucz$ where $z$ runs through some arbitrary transversal of $\Z_{F,+}^*/\Z_F^{*2}$.

Let $n = [F : \Q]$, $z \in \Z_{F,+}^*$, and $\xi \in I$. Then $\nrd \xi \in \nrd I = (ucz) \Z_F$, so $\alpha = (\nrd y)/(ucz) \in \Z_F$.
The arithmetic-geometric mean inequality implies
$$ n \le n (\N \alpha)^{1/n} \le \Tr \alpha = \varphi_{(ucz)^{-1}}(\xi). $$
Moreover, equality holds if and only if $1 = \N \alpha$ and $\alpha_v$ is independent of the real place $v$ of $F$, so equality holds if and only if  $\alpha = 1$.  Hence $\nrd(\xi) = ucz$ if and only if $\xi \in I$ satisfies $\varphi_{(ucz)^{-1}}(y) = n$ and is a shortest vector.
\end{proof}

To analyze the runtime of Algorithm \ref{findprincdef}, we first state and prove some preliminary results on enumeration of (short) vectors in lattices.

\begin{lemma} \label{LLLfixeddim}
There exists an algorithm which, given a lattice $L$ with basis vectors each of (square) norm at most $A \in \R_{>0}$, finds the shortest vectors in $L$ in deterministic time $O(\log^3 A)$ for fixed dimension $n$. 
\end{lemma}

\begin{proof}
This lemma is a consequence of the celebrated LLL-algorithm of Lenstra, Lenstra, and Lov\'asz \cite{LLL}: see e.g.\ Kannan \cite[Section 3]{Kannan}.
\end{proof}

\begin{rmk}
{\rm Using floating-point LLL, one could probably replace $\log^3 A$ with $\log^2 A$ in the lemma above, but this will have no impact on what follows so we neglect this possible improvement.}
\end{rmk}

We pause to prove two results which will be used in the next section.

\begin{lemma} \label{latticerandomelt}
A uniformly random element of a lattice $L$ with norm $\leq A$ can be computed in probabilistic probabilistic time (where the implied constant depends on the lattice $L$).
\end{lemma}

\begin{proof}
Choose a random vector $v \in \Z^n$ with each component $v_i$ chosen uniformly in the range $[-\sqrt{A},\sqrt{A}]$.  Using linear algebra over $\Q$, determine if $v \in L$; if so, output the vector, otherwise start again.  This test can be performed by first computing once (for the lattice $L$) the change of basis from $L$ to the standard basis and then multiplying the vector $v$ by this matrix.  The entries of the inverse will be of polynomial size for fixed dimension, and so each test can be performed in polynomial time.  The probability of success in each trial depends on the density of the lattice, given by its determinant.
\end{proof}

% A straightforward modification of exhaustive lattice enumeration using the Cholesky decomposition \cite[(2.8)]{FinckePohst}, \cite[\S 2.7.3]{Cohen} (instead choosing a random integer value in each successive position) yields the desired algorithm.  One performs $O(n^2)$ arithmetic operations, where $n=\dim L$, with integers of size $O(A)$, where the implied constant depends on the lattice $L$.  The result follows.

We apply this lemma to computing representatives of ideal classes, which we will use in the proof of Lemma \ref{findidealclass}.  

\begin{lemma} \label{randomidealinclass}
There exists an algorithm which, given $x \in \Z_{>0}$ sufficiently large and an ideal class $[\fraka] \in \Cl_{S_\infty}(\Z_F)$, computes a uniformly random ideal $\frakb \subset \Z_F$ with $[\frakb]=[\fraka]$ and $\N\frakb \leq x$, which runs in probabilistic polynomial time over a fixed totally real field $F$.
\end{lemma}

\begin{proof}
Let $\frakb \subset \Z_F$ be an ideal such that $N\frakb \leq x$ and $[\frakb] = [\fraka]$.  Then $\frakb\fraka^{-1}$ is principal.  By Lemma \ref{princZF}, there exists a generator $a$ for $\frakb\fraka^{-1}$ with 
\[ |a_v| \leq U \N(\frakb\fraka^{-1}) \leq U x / \N\fraka=T \] 
for all real places $v$ of $F$, where $U$ is an effectively computable constant depending only on $F$.  In particular, we have that
\[ \{\frakb : \N \frakb \leq x, [\frakb] = [\fraka] \} \subset \{ a\fraka^{-1} : a \in \Z_F, |a_v| \leq T\}. \]

Now each fractional ideal of $F$ is naturally a lattice under the Minkowski embedding $\sigma:F \to \R^n$ by $x \mapsto (x_v)_v$.  Denote by $\SV(L)$ the set of shortest vectors in a lattice $L$.  Then the map 
\begin{align*}
\{ a \in \Z_F: |a_v| \leq T\text{ for all $v$ and }a \in \SV(a\Z_F)\} &\to \{\frakb : \N \frakb \leq x\text{ and }[\frakb] = [\fraka] \} \\
a &\mapsto \frakb=a\fraka^{-1}
\end{align*}
is a surjective map of finite sets.  If we draw an element $a \in \Z_F$ on the left-hand side with probability $1/\#\SV(a\Z_F)$ and an ideal $\frakb$ on the right-hand side uniformly, then this map preserves probabilities.  

Therefore, we may compute a uniformly random ideal $\frakb$ as follows.  First, we compute a uniformly random element $a \in \fraka^{-1}$ such that $\N(a) \leq x/\N\fraka$ and let $\frakb=a\fraka$.  We find such an element $a$ by finding a uniformly random element $a \in \fraka^{-1}$ of (square) norm $\leq nT^2$ in the Minkowski embedding, which can be done in probabilistic polynomial time by Lemma \ref{latticerandomelt}.  If $\N a > T$, we return and compute another random $a$.  Given $a$ with $\N a \leq T$, we compute the shortest vectors in the lattice $a\Z_F$, which can be done in deterministic polynomial time by Lemma \ref{LLLfixeddim} (since $a$ is of polynomial size, the lattice $a\Z_F$ has a basis which is of polynomial size with the constant depending on the size of a basis for $\Z_F$).  If $a$ is not a shortest vector, we return and compute another element $a$; if $a$ is a shortest vector, then we keep $a$ with probability $1/\#\SV(a\Z_F)$ and otherwise return.
\end{proof}

We conclude this section by analyzing the runtime of the principalization algorithm.  

\begin{proposition} \label{dembdonn}
Algorithm \textup{\ref{findprincdef}} runs in deterministic polynomial time in the size of the input over a fixed totally real field $F$.
\end{proposition}

\begin{proof}
In Step 1, we must solve Problem  (\textsf{ClassUnitGroup}) for the field $F$ which can be done in constant time.  In Step 2, we note that $\log C$ is of polynomial size in the input, with the constant depending only on $F$.  In Step 3, there are at most $2^{n-1}$ totally positive units $z$.  For each, in Step 3(a) the $\Z$-lattice $aI$ has basis vectors with norm of polynomial size in the input pseudobasis for $I$, with the constant depending on $F$.  In Step 3(b) we find the shortest vectors in this lattice, and by Lemma \ref{LLLfixeddim} this can be performed in deterministic polynomial time in fixed dimension.
\end{proof}

\begin{rmk}
{\rm We now present an alternative approach to Corollary \ref{isconjbyconnectingideal} which solves Problem \ref{conjugacytest} \textup{(\textsf{IsConjugate})} for a definite algebra by constructing an isomorphism directly. Let 
\[ O(B) = \{ \varphi \in \End_F(B) : \nrd(\varphi(x)) = \nrd(x) \text{ for all } x \in B \} \] 
where $\End_F(B)$ denotes the endomorphisms of $B$ as an $F$-vector space, and 
\[ SO(B) = \{ \varphi \in O(B) : \det(\varphi) = 1\}. \] 
It is well known (see e.g.\ Dieudonn\'e \cite[Appendix IV, Proposition 3]{Dieudonne}) that
\[ SO(B) = \{ \phi : x \mapsto \nu\delta x \delta^{-1} : \nu,\delta \in B^*,\ \nrd(\nu)=1 \}. \]
Since conjugation is an order-preserving anti-automorphism of determinant $-1$, we have the following result: Two orders $\calO$ and $\calO'$ of $B$ are conjugate if and only if there exists $\varphi \in O(B)$ such that $\varphi(\calO) = \calO'$.

Therefore, such an isometry can be found using the approach described by Plesken and Souvignier \cite{Souvignier}.}
%  This solves Problem \ref{conjugacytest} (\textsf{IsConjugate}) in the totally definite case.
\end{rmk}

\begin{rmk}
{\rm For a definite quaternion order $\calO$, we define its \emph{theta series} simply to be the theta series of the corresponding lattice under the quadratic form $\Tr \nrd$, i.e., 
\[ \theta(\calO;q)=\sum_{\xi \in \calO} q^{\Tr\nrd(\xi)} \in \Z[[q]]. \]
Isomorphic quaternion orders have the same theta series, and so given a complete set of representatives $\{\calO_i\}$ for the set of conjugacy classes of (definite) Eichler orders of a given level, one can identify the conjugacy class of a given order $\calO$ by comparing the first few coefficients of their theta series.  This approach is more efficient in practice to show that two orders are not conjugate than to test for conjugacy directly.

In this way, one can also show that two right ideals of an order $\calO$ are not isomorphic, by showing that their left orders are not conjugate, and this idea can be used in practice to speed up the enumeration of ideal classes.}
\end{rmk}

\section{Computing the class set for definite orders}

In this section, we discuss Problem \ref{computeidealclasses} (\textsf{ClassSet}) in the totally definite case.  The key algorithm is the $S$-neighbors algorithm (Algorithm \ref{Sneighbors}), but see also Algorithm \ref{combined}, and its running time is analyzed in Proposition \ref{conjruntime}.  The main result (the definite case of Theorem B) appears as Theorem \ref{thmb}.  As in the previous section, let $B$ be a totally definite quaternion algebra of discriminant $\frakD$ and let $\calO \subset B$ be an Eichler order of level $\frakN$.  

One way to solve Problem (\textsf{ClassSet}) is to simply enumerate the right invertible fractional $\calO$-ideals in some way, building a list of ideal classes by testing each ideal (an instance of Problem (\textsf{IsIsomorphic})) to see if it is isomorphic to any ideal in the list, and stopping when one knows that the list is complete.  Ultimately, our methods to solve Problem (\textsf{ClassSet}) will be a variant of just this simple idea.  We use the mass formula (Proposition \ref{massformula}) and the accompanying solution to Problem (\textsf{ClassNumber}) and stop when equality holds.  (See also Remark \ref{usealtmassformula} below for the use of the alternate mass formula.)

\begin{rmk}
{\rm In the commutative case, one bounds the norm of the ideals in a generating set for $\Cl \Z_F$ using a Minkowski-like bound (or the Bach bound, assuming the GRH).  One can in a similar way use the geometry of numbers in the case of quaternion algebras: the first author \cite{Kirschmer} has proven that if $\calO$ is maximal, every ideal class in $\Cl \calO$ is represented by an ideal $I$ with 
\[ \nrd(I) \leq M_n d_F^2 \sqrt{\N \frakD} \]
where $M_n$ is a constant growing exponentially with $n$.  Although this bound is effective, just as in the commutative case, it would be all but useless in practice except in the simplest cases.}
\end{rmk}

To enumerate ideals, we note that the group $\Cl \Z_F$ acts on the set of isomorphism classes of fractional right $\calO$-ideals by multiplication (on the left or right).  We say that an invertible right (integral) $\calO$-ideal $I$ is \emph{primitive} if $I$ is not contained in any nontrivial two-sided ideal of the form $J=\fraka \calO$ with $\fraka \subset \Z_F$.  We note that $I \subset J$ if and only if $IJ^{-1} \subset \calO$ is integral, and if $I \subset J$ then $\nrd(J) \mid \nrd(I)$ so $I$ is contained in only finitely many two-sided ideals $J$.  

Therefore, to enumerate all right $\calO$-ideals, we enumerate the products $JI$ of primitive right $\calO$-ideals $I$ and two-sided ideals $J=\fraka \calO$.  To enumerate primitive right ideals, we employ the following lemma.  For a commutative ring $R$, we denote by
\[ \PP^1(R)=\{(x,y) : xR+yR=R\}/R^* \]
the points of the projective line over $R$, and denote by $(x:y)$ the class of $(x,y)$ in $\PP^1(R)$.

\begin{lemma} \label{idealsnorma}
Let $\fraka$ be an ideal of $\Z_F$ coprime to $\frakD\frakN$.  Then the set of primitive right invertible $\calO$-ideals of norm $\fraka$ are in bijection with $\PP^1(\Z_F/\fraka)$.  Explicitly, given a splitting
\[ \phi_\fraka:\calO \hookrightarrow \calO \otimes_{\Z_F} \Z_{F,\fraka} \cong M_2(\Z_{F,\fraka}) \to M_2(\Z_F/\fraka) \]
the bijection is given by
\begin{align*}
\PP^1(\Z_F/\fraka) &\to \{ I \subset \calO : \nrd(I)=\fraka \} \\
(x:y) &\mapsto I_{(x:y)}=\phi_\fraka^{-1}\begin{pmatrix} x & 0 \\ y & 0 \end{pmatrix}\calO + \fraka\calO.
\end{align*}
\end{lemma}

\begin{proof}
We may refer to Vign\'eras \cite[Chap.\ II, Th\'eor\`eme 2.3(3)]{Vigneras}); for convenience, we give a proof here.

It is clear that the ideal $I_{(x:y)}$ is a right $\calO$-ideal of norm $\fraka$; it is invertible since it is locally principal, and it is primitive since the only two-sided ideals with norm dividing $\fraka$ are of the form $\frakp \calO$ with $\frakp \mid \fraka$, and $\begin{pmatrix} x & 0 \\ y & 0 \end{pmatrix} \not\in \frakp M_2(\Z_{F,\fraka})$ for any such $\frakp$.

To prove the lemma, it suffices to construct an inverse map.  Let $I$ be a right invertible $\calO$-ideal of norm $\fraka$.  Then $\phi_\fraka(I)$ is a right principal $M_2(\Z_F/\fraka)$-ideal, say, generated by $\xi \in M_2(\Z_F/\fraka)$.  We have $\det(\xi)=0$, so the kernel 
\[ V=\{v \in (\Z_F/\fraka)^2 : v\xi = 0\} \] 
of $\xi$ acting on the right is nonzero.  We claim in fact that $V$ is a free $\Z_F/\fraka$-module of rank $1$.  Indeed, for each prime $\frakp \mid \fraka$, we have $V/\frakp V$ one-dimensional since $I$ is primitive and hence $\xi \not\equiv 0 \pmod{\frakp}$; it follows from a Hensel lift that $V/\frakp^e V$ is also one-dimensional, and the claim then follows by the Chinese remainder theorem.  This argument also shows that a generator $v=(x,y)$ of $V$ has $(x:y) \in \PP^1(\Z_F/\fraka)$ (and this element is unique), since this is true modulo $\frakp$ for all $\frakp \mid \fraka$.  Thus the association $I \mapsto (-y:x)$ is well-defined, and it is easy to see that this furnishes the desired inverse.
\end{proof}

One natural way to enumerate primitive ideals $I$ would be to order them by $\N(\nrd(I))$, the absolute norm of their reduced norm.  Alternatively, one may restrict the set of possible reduced norms as follows.

\begin{proposition}[Strong approximation] \label{strongapprox}
Let $S$ be a nonempty set of (finite) prime ideals of $\Z_F$ coprime to $\frakD$ which generate $\Cl_{S_\infty} \Z_F$.  Then there exists a set of representatives $I$ for $\Cl \calO$ such that $\nrd(I)$ is supported in $S$.
\end{proposition}

\begin{proof}
See Vign\'eras \cite[Theor\`eme III.4.3]{Vigneras} and the accompanying discussion.
\end{proof}

Let $\fraka$ be a squarefree ideal of $\Z_F$, and let $I,J$ be right invertible $\calO$-ideals.  Then $J$ is said to be a \emph{$\fraka$-neighbor} of $I$ if $I \supset J$ and $\nrd(J)=\fraka\nrd(I)$, or equivalently if $[I:J]=\fraka^2$, the index taken as $\Z_F$-lattices.  Following Schulze-Pillot \cite{Schulze-Pillot}, Pizer \cite{Pizer}, Kohel \cite{Kohel}, and Mestre \cite{Mestre}), we enumerate primitive ideals by iteratively enumerate the $\frakp$-neighbors as follows.  

\begin{alg}[$S$-neighbors] \label{Sneighbors}
Let $\calO$ be an Eichler order of level $\frakN$ in a quaternion algebra $B$ of discriminant $\frakD$, and let $S$ be a nonempty finite set of prime ideals of $\Z_F$ coprime to $\frakD\frakN$ that generate $\Cl_{S_\infty} \Z_F$.  This algorithm solves Problem \textup{(\textsf{ClassSet})}.
{\rm \begin{enumerate}
\item Solve \textup{(\textsf{ClassNumber})} as in Theorem \ref{existstotdef} and let $H=\#\Cl \calO$.
\item Compute $\Cl \Z_F$ and compute a set $C$ of representatives of the set of ideals of $\Z_F$ supported at $S$ modulo principal ideals.
\item Initialize $\calI := \emptyset$, $\calF := \{\calO\}$, and $\calF_{\text{new}} := \emptyset$.
\item For each $I \in \calF$, compute the primitive $\frakp$-neighbors $I'$ of $I$ using Lemma \ref{idealsnorma} for $\frakp \in S$.  For each such $I'$, determine if $I'$ is isomorphic to any ideal in $\calI \cup \calF_{\text{new}}$ by an algorithm to solve \textup{(\textsf{IsIsomorphic})}; if not, add $JI'$ to $\calF_{\text{new}}$ for all $J \in C$.
\item Set $\calI := \calI \cup \calF$, $\calF := \calF_{\text{new}}$ and $\calF_{\text{new}} := \emptyset$.  If $\#\calI=H$, return $\calI$, otherwise return to Step 4.
\end{enumerate}}
\end{alg}

\begin{proof}[Proof of correctness]
The algorithm enumerates all right invertible $\calO$-ideals with norm supported at $S$ by the discussion preceding Lemma \ref{idealsnorma} and so will find all ideal classes by Proposition \ref{strongapprox}.
\end{proof}

Our algorithm was inspired by the implementation due to Kohel in the computer algebra system \textsf{Magma} \cite{Magma} for definite quaternion algebras over $\Q$.  

\begin{rmk} \label{usealtmassformula}
{\rm In Algorithm \ref{Sneighbors}, one can alternatively use the mass formula as in Remark \ref{altmassformula}; here, one trades the difficulty of computing the class number directly with computing the unit index $[\calO^*:\Z_F^*]$ for each order $\calO=\calO_L(I)$.  To compute this unit index, we consult Vign\'eras \cite[Theorem 6]{VignerasSimpl}.  If all totally positive units of $\Z_F$ are squares, then the group $\calO_L(I)^*/\Z_F^*$ is just the set of shortest vectors in $\calO_L(I)$, i.e., the reduced norm $1$ subgroup $\calO_1^*$ modulo $\{\pm 1\}$, and this group can be computed by lattice enumeration.  Otherwise, if the norm $1$ subgroup is not cyclic, then $\calO^*/\Z_F^*$ is at most an extension of
$\calO_1^*/\{\pm 1\}$ of index $2$ and we can explicitly write down candidates for this extension.  If the norm $1$ subgroup is cyclic, then we must fall back on a lattice search in $\calO$.  }
\end{rmk}

To analyze the running time of Algorithm \ref{Sneighbors}
, we now examine the distribution of ideal classes among ideals following this enumerative strategy.  Let $S$ be a nonempty finite set of ideals of $\Z_F$ coprime to $\frakD\frakN$ (as in Algorithm \ref{Sneighbors}
).  Let $G(S)$ be the graph with vertices $[I] \in \Cl(\calO)$ and edges as follows: for a vertex $[I]$, choose a representative ideal $I$, and for each primitive $\frakn$-neighbor $J$ with $\frakn \in S$, draw an edge from $[I]$ to $[J]$.  By Lemma \ref{idealsnorma}
, $\Gamma(S)$ is $k$-regular (the out degree of each vertex is $k$), where $k=\sum_{\frakn \in S} \Phi(\frakn)$ and 
\[ \Phi(\frakn)=\prod_{\frakp^e \parallel \frakn} N\frakp^{e-1} (N\frakp+1). \]
When $S=\{\frakn\}$ consists of a single ideal, we abbreviate $G(S)$ by simply $G(\frakn)$.  

Suppose now that the map $\nrd : S \to \Cl_{S_\infty}(\Z_F)$ is surjective (i.e., that norms of ideals in $S$ cover all narrow ideal classes of $\Z_F$).  Then the set of vertices of the quotient $G(S)$ is a set of representatives for $\Cl \calO$ by Proposition \ref{strongapprox}
.

The graphs $G(\frakp)$ are known in many cases to be \emph{Ramanujan graphs}, graphs whose eigenvalues other than $\pm k$ have absolute value at most $2\sqrt{k-1}$, and are therefore a type of \emph{expander graphs}.  In the simplest case where $F=\Q$ and $S=\{p\}$, they were first studied by Ihara, then studied in specific detail by Lubotzky, Phillips, and Sarnak \cite{LPS} and Margulis \cite{Margulis}: the nontrivial spectrum of $G(p)$ can be identified with the spectrum of the Hecke operators acting on the space of cusp forms of weight 2 on $\Gamma_0(p)$, and the eigenvalue bound then follows from the Ramanujan-Petersson conjecture, a consequence of the Eichler-Shimura isomorphism and Deligne's proof of the Weil conjectures (though strictly speaking, one only needs the result in weight $2$, so could be derived from earlier results of Eichler, Igusa, Weil, and Shimura).   Charles, Goren, and Lauter \cite{CGL} have shown that the same is true for graphs $G(\frakl)$ under the following set of hypotheses: $F$ is a totally real field of narrow class number one, $B \cong B_p \otimes_\Q F$ where $B_p$ is the quaternion algebra of discriminant $p$ over $\Q$ and $p$ is unramified in $F$, and $\frakl$ is a prime ideal of $\Z_F$ coprime to $p$.  We note that there are further constructions due to Jordan and Livn\'e \cite{JL}.  

Let $T(S)$ be the adjacency matrix of $G(S)$.  We begin by discussing the interpretation of $T(S)$ as representing the action of a Hecke operator on the space of Hilbert modular forms.   We refer the reader to work of Demb\'el\'e and the second author \cite{DemVoi} for a computational point of view on the statements below, and the references therein for further information.  

The graph $G(\frakn)$, for an ideal $\frakn$ coprime to $\frakD\frakN$, represents the action of a Hecke operator on a certain space of quaternionic modular forms of level $\frakN$ for the quaternion algebra $B$.  By the Jacquet-Langlands correspondence \cite{JacquetLanglands} (see also e.g.\ work of Hida \cite{HidaCM}), this space as a Hecke module (finite-dimensional $\C$-vector space equipped with the action of Hecke operators) is isomorphic to a subspace $M_2$ of the space of Hilbert modular forms of parallel weight $2$ and level $\frakD\frakN$.  Therefore, we may identify eigenvectors of $T(S)$ with Hilbert modular eigenforms.  

Since the adjacency matrix $T(S)$ of the graph $G(S)$ is the sum of the adjacency matrices of the graphs $G(\frakn)$ for $\frakn \in S$, the matrix $T(S)$ represents the action of the operator
\[ T_S = \sum_{\frakn \in S} T_\frakn. \]
For $\frakn$ coprime to $\frakD\frakN$, the operators $T_\frakn$ are semisimple and pairwise commute, and so there is a set of common eigenforms for all $T_\frakn$ which span $M_2$.  The space $M_2=S_2 \oplus E_2$ breaks up accordingly into a space of Eisenstein series $E_2$ and a space of cusp forms $S_2$.

We now apply eigenvalue bounds coming from geometry.  We first discuss the Eisenstein series.  
The space of Eisenstein series $E_2$ corresponds, under Jacquet-Langlands, to the standard Eisenstein series and its twists by Hilbert class characters, and is spanned by those vectors which ``factor through the reduced norm'', i.e., vectors of the form $a_I [I]$ with $a_I=a_J$ whenever $[\nrd(I)]=[\nrd(J)] \in \Cl_{S_\infty}(\Z_F)$.  The eigenvectors are then given concretely as follows: for every irreducible character $\chi$ of $\Cl_{S_\infty}(\Z_F)$, the vector $v_\chi=( \chi([\nrd(I)]) )_{[I] \in \Cl(O)}$ is an eigenvector of $T(\frakp)$ with eigenvalue 
\[ \chi([\frakp])\Phi(\frakp) \] 
and thus an eigenvector of $T(S)$ with eigenvalue
\[ \sum_{\frakn \in S} \chi([\frakn])\Phi(\frakn). \] 

Now we turn to eigenvalue bounds for the remaining space $S_2$ of cusp forms.  (We obtain exactly those cusp forms of level $\frakD\frakN$ which are new at all primes dividing $\frakD$, but this is not relevant for our bounds.)  The general result we will need is due to Livn\'e \cite{Livne}, as follows.

\begin{theorem}  \label{ramanujan}
Let $\frakp$ be a prime with $\frakp \nmid d_F \frakD\frakN$.  Then the eigenvalues of $T_\frakp$ acting on the space of cusp forms $S_2$ are bounded by $2\sqrt{N\frakp}$.  In particular, the graph $G(\frakp)$ is Ramanujan if $\frakp \nmid d_F \frakD\frakN$.
\end{theorem}

\begin{proof}
See the proof of the Ramanujan-Petersson conjecture for Hilbert modular forms by Livn\'e \cite[Theorem 0.1]{Livne} and the accompanying discussion.
\end{proof}

Suppose from now on that $\frakn$ is squarefree.  Since $T_\frakn = \prod_{\frakp \mid \frakn} T_\frakp$, it follows that the eigenvalues of $T_\frakn$ acting on $S_2$ are bounded by $\prod_{\frakp \mid \frakn} 2\sqrt{N\frakp}$ and accordingly the eigenvalues of $T_S$ on $S_2$ are bounded by
\[ \sum_{\frakn \in S} \prod_{\frakp \mid \frakn} 2\sqrt{N\frakp}. \]

Finally, the space $M_2=E_2 \oplus S_2$ is equipped with a natural inner product coming from the Jacquet-Langlands correspondence called the Petersson inner product, under which cuspidal eigenforms are mutually orthogonal and $E_2$ is the orthogonal complement of $S_2$.  This inner product on the free $\Z$-module spanned by $[I] \in \Cl(\calO)$ is defined by 
\[ \la [I], [J] \ra = \delta_{[I],[J]} w_I/2 \]
where $w_I=\# \calO_{L}(I)^\times/\Z_F^\times$ and $\delta_{[I],[J]}=0$ or $1$ according as $[I] \neq [J]$ or $[I]=[J]$.

More is true: in fact, the matrix $T(S)$ is normal with respect to this inner product.  Since this is a key point in our proof and the source of our original mistake, we give a proof.

\begin{claim*}
$T=T(S)$ is normal with respect to the Petersson inner product.
\end{claim*}

\begin{proof}
Let $I_1,\dots,I_H$ be representatives of the right ideal classes of $\calO$; for $i=1,\dots,H$ let $\calO_i = \calO_L(I_i)$, $\Gamma_i=\calO_i^\times/\Z_F^\times$ and $w_i=\#\Gamma_i$.  We have shown that $T$ has a basis of eigenvectors; we show that these eigenvectors are orthogonal with respect to the Petersson inner product: this is true already for the ones corresponding to cusp forms, so it is enough to show that the vectors $v_\chi$ are orthogonal.  

Since the eigenvalue of $T(\frakp)$ on $v_\chi$ is $\chi([\frakp])\Phi(\frakp)$, by the Chebotarev density theorem we can choose primes $\frakp$ so that the eigenvalues of $v_\chi$ are distinct for each $\chi$.  Since there is a common basis of eigenvectors for all $T(\frakp)$, we conclude that this common eigenbasis is unique (unordered).  Therefore, if we show that $T(\frakp)$ is normal for each prime $\frakp$, so the eigenbasis is orthogonal, then this is also an eigenbasis for $T(S)$, so $T(S)$ is normal.

So suppose that $S=\{\frakp\}$.  Let $m$ be the order of $[\frakp] \in \Cl_{S_\infty}(\Z_F)$.  We will show that $T^m$ is self-adjoint.  

An inclusion $J=xI_i \subset I_j$ with $\nrd(J)=\nrd(I_j)\frakp$ corresponds to an element $x \in I_jI_i^{-1}$ with 
$\nrd(xI_i I_j^{-1})=\frakp$ and $x$ is uniquely determined by $I_i$ and $J=x I_i$ up to multiplication on the right by an element of $\calO_i^\times$.  Therefore by definition, the coefficient $T(S)_{i,j}$ is equal to the cardinality of the set
\[ \{x \in I_j I_i^{-1} : \nrd(x I_i I_j^{-1}) = \frakp \} / \calO_i^\times \]
which is 
\[ \frac{1}{w_i} \#\{x \in I_j I_i^{-1} : \nrd(x I_i I_j^{-1}) = \frakp \} / \Z_F^\times. \]

We similarly see that the coefficient $w_i (T^m)_{i,j}$, interpreted as the (weighted) number of paths in $G(S)$ of length $m$, is equal to
\begin{align*}
w_i (T^m)_{i,j} &= \# \Theta_{i,j}(\frakp^m) \\
&=
\# \{x \in I_j I_i^{-1} : \nrd(x I_i I_j^{-1}) = \frakp^m\} / \Z_F^\times. 
\end{align*}
If $(T^m)_{i,j}=0$, then $(T^m)_{j,i}=0$ as well by the argument below.  So suppose that $\Theta_{i,j}(\frakp^m)$ is nonempty.  Then $[\nrd(I_i I_j^{-1})]=[\frakp^m \nrd(x)^{-1}]=1$, so there exists a totally positive $a_{ij} \in \Z_F$ such that $a_{ij}\Z_F = \nrd(I_i I_j^{-1})$.  
We define a map
\begin{align*}
\Theta_{i,j}(\frakp^m) &\to \Theta_{j,i}(\frakp^m) \\
x &\mapsto a_{ij}\overline{x}
\end{align*}
where the bar $\overline{\phantom{x}}$ indicates conjugation in $B$.  Indeed, if $x \in \Theta_{i,j}(\frakp^m)$ then $x \in I_jI_i^{-1}$ so 
\[ \overline{x} \in \overline{I_i}^{-1} \overline{I_j} = I_i I_j^{-1} \nrd(I_jI_i^{-1}) \]
so $a_{ij}\overline{x} \in I_i I_j^{-1}$; and
\[ \nrd(a_{ij}\overline{x} I_j I_i^{-1})=\nrd(x)a_{ij}^2\nrd(I_j I_i^{-1}) = \nrd(x I_i I_j^{-1}) =
\frakp^m. \]
Therefore, this map is a bijection, so $w_i (T^m)_{i,j} = w_j (T^m)_{j,i}$, and thus $T^m$ is self-adjoint as claimed.

Now, in view of our comments above, the (uniquely defined) eigenbasis for $T$ is an eigenbasis for $T^m$ and so is pairwise orthogonal by the spectral theorem, thus $T$ is normal.
\end{proof}

Define the \emph{distance} $d(v,w)$ between two vertices $v,w$ in the graph $G(S)$ to be the length of the shortest path between them, and define the \emph{diameter} of $G$ to be $D(G)=\max_{v,w} d(v,w)$.  To estimate the run time of the $S$-neighbors algorithm (\ref{Sneighbors}
), we bound the diameter in a special case.  Recall that for coprime ideals $\frakD,\frakN$, by Remark 5.2 we have the mass formula
\[ M(\frakD,\frakN)=\sum_{[I] \in \Cl(\calO)} \frac{1}{w_I}. \]
Let $w_{\textup{max}}=\max_{[I] \in \Cl(\calO)} w_I$.

\begin{proposition} \label{conjruntime}
Let $S$ be a finite set of primes $\frakp$ of $\Z_F$ representing $\Cl_{S_\infty}(\Z_F)$ with $\frakp \nmid d_F\frakD\frakN$.   Let $\frakp$ be the smallest prime in $S$.  Let $m$ be the smallest positive integer which is congruent to $1$ modulo the exponent of $\Cl_{S_\infty}(\Z_F)$ which satisfies
\[ m > \frac{\log (M(\frakD, \frakN)w_{\textup{max}}/h^+(\Z_F))}{\log \bigl((N\frakp+1)/(2\sqrt{N\frakp})\bigr)}. \]
Then $D(G(S)) \leq m$.
\end{proposition}

\begin{proof}
We follow the proof of Chung \cite{Chung}.  We want to find the minimum value of $m$ such that $T(S)^m$ has all entries nonzero.

As in the proof of the claim, let $I_1,\dots,I_H$ be representatives of the right ideal classes of $\calO$; for $i=1,\dots,H$, let $w_i=\#(\calO_L(I_i)/\Z_F^\times)$ and $\frakp_i=\nrd(I_i)$, and let $W$ be the diagonal matrix with entries $1/w_i$.  Abbreviate $M=M(\frakD,\frakN)$.

Let $1 \leq r,s \leq H$ be indices.  We have $T(S)=\sum_{\frakp \in S} T(\frakp)$ so 
\[ (T(S)^m)_{r,s} \geq \sum_{\frakp \in S} (T(\frakp)^m)_{r,s}. \]
For indices $r,s$, let $\frakp_{r,s}$ be such that $[\frakp_{r,s}]=[\frakp_s^{-1} \frakp_r]$, possible since $\frakp$ covers $\Cl_{S_\infty}(\Z_F)$.  Then $(T(S)^m)_{r,s} \geq (T(\frakp)^m)_{r,s}$.  So without loss of generality, we may assume $S=\{\frakp_{r,s}\}$.  To simplify, abbreviate $\frakp=\frakp_{r,s}$ and $T=T(S)$.

Let $\chi_1,\dots,\chi_{h^+}$ be the irreducible characters of $\Cl_{S_\infty}(\Z_F)$, so $h^+=h^+(\Z_F)$.  For $i=1,\dots,h^+$, let 
\[ u_i = \frac{1}{\sqrt{M}}( \chi_i ([\frakp_j]) )_{j=1,\dots,H} \]
be the eigenvector corresponding to an Eisenstein series, with eigenvalue 
\[ \chi_i([\frakp])(N\frakp+1). \] 
Then the eigenvectors $u_i$ are orthonormal with respect to the matrix $W$, i.e., $u_i W u_j^* = \delta_{ij}$ for $i,j=1,\dots,h^+$, where ${}^*$ denotes conjugate transpose.

Complete $u_1,\dots,u_{h^+}$ to a basis $u_1,\dots,u_H$ of (row) eigenvectors for $T$ normalized so that $u_i W u_i^* = 1$ for all $i$.  Then by the orthogonality relations of the Petersson inner product, we have $u_i W u_j^* = 0$ for all $i \neq j$.  Let $U$ be the matrix with rows $u_i$.  Then $UWU^*=1$, and since $W$ is diagonal with positive diagonal entries, if $V=U\sqrt{W}$ (taking positive square roots), then $V^{-1}=V^*$ is unitary.  

Now let $D$ be the diagonal matrix whose diagonal entries are the eigenvalues $\lambda_1,\dots,\lambda_H$ of $u_1,\dots,u_H$.  
Let $T'=\sqrt{W}^{-1} T \sqrt{W}$.  Then for any $m \geq 0$, since $W$ is diagonal we have
\[ (T')^m = \sqrt{W}^{-1} T^m \sqrt{W} = \sqrt{W}^{-1} U^{-1} D^m U \sqrt{W} = V^* D^m V. \] 
Note that $(T')^m$ has nonnegative real entries, and since $T^m$ has positive entries if and only if $(T')^m$ does, we may proceed using the argument in Chung's proof.  

Let $v_i$ be the rows of $V$.  We have $T'=\sum_{i=1}^{H} \lambda_i v_i^* v_i$ and so for $r,s \in 1,\dots,H$, the $r,s$-entry of $(T')^m$ is equal to 
\[ ((T')^m)_{r,s} \geq \sum_{i=1}^{H} \lambda_i^m (v_i^* v_i)_{r,s} = 
\sum_{i=1}^{h^+} \lambda_i^m (v_i^* v_i)_{r,s} + \sum_{i=h+1}^{H} \lambda_i^m (v_i^* v_i)_{r,s}. \]
Then the first term is equal to
\begin{align*}
\sum_{i=1}^{h^+} \lambda_i^m (v_i^* v_i)_{r,s} &= \frac{1}{M\sqrt{w_rw_s}}\sum_{i=1}^{h^+}
\overline{\chi}_i([\frakp_r]){\chi}_i([\frakp_s])\bigl(\chi_i([\frakp])(N\frakp+1)\bigr)^m \\
&=\frac{1}{M\sqrt{w_rw_s}} (N\frakp +1)^m \sum_{i=1}^{h^+} \chi_i([\frakp_s \frakp^m])\overline{\chi}_i([\frakp_r]).
\end{align*}
Now the theory of character sums says that for a finite abelian group $G$ and $x,y \in G$, the sum over the set of all irreducible characters $\widehat{G}$
\[ \sum_{\chi \in \widehat{G}} \chi(x)\overline{\chi}(y) = \begin{cases} \#G, & \text{if $x=y$}; \\ 0, & \text{if $x \neq y$}. \end{cases} \]
Therefore the inner first sum is $h$ if $[\frakp_s\frakp^m] = [\frakp_r]$ and $0$ otherwise.  

We now suppose $m \equiv 1 \pmod{\exp(\Cl_{S_\infty}(\Z_F))}$ where $\exp$ denotes the exponent of the group.  Recall that we chose $\frakp$ so that $[\frakp]=[\frakp_s^{-1}\frakp_r]$.  Therefore we have
\begin{align*}
\sum_{i=1}^{h^+} \lambda_i (v_i^* v_i)_{r,s}  \geq \frac{h^+}{M\sqrt{w_rw_s}} (N\frakp+1)^m.
\end{align*}
Let $\lambda=2\sqrt{N\frakp}$.  Then $|\lambda_i| \leq \lambda$ for $i=h+1,\dots,H$.  Continuing with the second sum, we reason as Chung that 
\begin{align*} 
\left|\sum_{i={h^+}+1}^{H} \lambda_i^m (v_i^* v_i)_{r,s}\right| &\leq \lambda^m \left(1 - \sum_{i=1}^{h^+} |(v_i)_r|^2 \right)^{1/2}
\left(1 - \sum_{i=1}^{h^+} |(v_i)_s|^2 \right)^{1/2} \\
&= \lambda^m \left(1 - \frac{h^+}{Mw_r} \right)^{1/2}
\left(1 - \frac{h^+}{Mw_s} \right)^{1/2} \leq \lambda^m.
%%&\leq \lambda^m\left(1 + \frac{h^2}{w^2w_rw_s}\right)^{1/2}.
\end{align*}
Let $t = h^+/(M\sqrt{w_rw_s})$ and $k = N\frakp+1$.  Then putting these together, we have
\[
((T')^m)_{r,s} \geq t k^m - \lambda^m. \]
Thus $((T')^m)_{r,s} > 0$ whenever 
\[ m > \frac{\log(t)}{\log (\lambda/k)}=\frac{\log(M\sqrt{w_r w_s}/h^+)}{\log(k/\lambda)}. \]

To obtain an estimate for all indices $r,s$, let $w_{\textup{max}}=\max_i w_i$, let $\frakp_{\textup{min}}=\min_{\frakp \in S} N\frakp$ and $k_{\textup{min}}=N\frakp_{\textup{min}}+1$, and $\lambda_{\textup{min}}=2\sqrt{N\frakp_{\textup{min}}}$.  Then for
\[ m > \frac{\log (Mw_{\textup{max}}/h^+)}{\log k_{\textup{min}}/ \lambda_{\textup{min}}} \]
with $m \equiv 1 \pmod{\exp(\Cl_{S_\infty}(\Z_F))}$  we have that every coefficient of $(T')^m$ is positive.  This completes the proof of the proposition.
\end{proof}

Before proving the definite case of Theorem B, we prove one further lemma.

\begin{lemma} \label{findidealclass}
There exists an algorithm which, given an ideal class $[\fraka] \in \Cl_{S_\infty}(\Z_F)$ and a finite set $T$ of primes of $\Z_F$, computes in probabilistic polynomial time a prime $\frakp \not\in T$ such that $[\frakp]=[\fraka]$ and $\N\frakp=O( \log^2 \prod_{\frakq \in T} \N\frakq)$ for a fixed totally real field $F$.
\end{lemma}

\begin{proof}
We use a (weak) version of the effective Chebotarev density theorem due to Lagarias and Odlyzko, applied to the strict Hilbert class field of $F$.  For an ideal class $[\fraka] \in \Cl_{S_\infty}(\Z_F)$, define the counting function 
\[ \pi_{[\fraka]}(x) = \#\{\frakp \subset \Z_F \text{ prime} : [\frakp]=[\fraka]\text{ and }\N\frakp \leq x\} \]
for $x \in \R_{\geq 2}$.  Then, by Lagarias and Odlyzko \cite[Theorems 1.3--1.4]{LagariasOdlyzko}, there exists an effectively computable constant $x_0>0$ which depends only on $F$ such that for all $[\fraka] \in \Cl_{S_\infty}(\Z_F)$ and all $x \geq x_0$, we have
\begin{equation} \label{effcheb}
\left|\pi_{[\fraka]}(x) - \frac{1}{h}\frac{x}{\log x}\right| \leq \frac{1}{2h}\frac{x}{\log x}
\end{equation}
where $h=\#\Cl_{S_\infty} \Z_F$.  
(In fact, they prove a certainly stronger bound on the error and show that $\log x_0=(d_F n)^{O(1)}$, but this weak version suffices for our purposes.)  

We begin with two precomputation steps.  We first compute the constant $x_0$, the strict class group $\Cl_{S_\infty} \Z_F$, and factor the discriminant $d_F$.  Next, for each (rational) prime $p \leq \Z_F$ with $p \leq x$, we factor $p\Z_F$ and see if there exists a prime $\frakp \mid p$ such that $[\frakp]=[\fraka]$ and $\frakp \not\in T$.  If so, we return the ideal $\frakp$.  If not, then for all primes $\frakp$ with $[\frakp]=[\fraka]$ we have $\frakp \in T$.  

Let $\N(T)=\prod_{\frakq \in T} \N\frakq$ and let $x = \max(x_0, \log^2 \N(T), (4h)^4)$.  Clearly $x = O(\log^2 \N(T))$.  Then 
\[ \#T \leq \N(T) \leq \sqrt{x} \leq \frac{x}{4h \log x}. \]  
By (\ref{effcheb}), we have
\[ \frac{\pi_{[\fraka]}(x)- \#T}{x} > \frac{1}{2h \log x} - \frac{1}{4h\log x} = \frac{1}{4h\log x}. \]
We have that $\#\{\frakb \in [\fraka] : \N\frakb \leq x\} \leq cx$ for some constant $c>0$.  Therefore, by employing the algorithm of Lemma \ref{randomidealinclass}
, we can find in probabilistic polynomial time a prime ideal $\frakp$ with $[\frakp]=[\fraka]$, $\N\frakp \leq x$, and $\frakp \not\in T$.
\end{proof}

We now prove the definite case of Theorem B.

\begin{theorem} \label{thmb}
There exists an algorithm to solve Problem \textup{(\textsf{ClassSet})} for definite orders with factored discriminant over a fixed field $F$ which runs in probabilistic polynomial time in the output size.
\end{theorem}

\begin{proof}
We use Algorithm \ref{Sneighbors}
.  We must, for a given order $\calO$ of discriminant $\frakd=\frakD\frakN$, choose a set $S$ of primes $\frakp$ of $\Z_F$ which represents the elements of the narrow class group $\Cl_{S_\infty} \Z_F$ such that $\frakp \nmid d_F\frakD\frakN$.  We do this by first computing a set of representatives for $\Cl_{S_\infty} \Z_F$ (which runs in constant time for fixed $F$) and then using the algorithm of Lemma \ref{findidealclass} with $T=\{\frakp : \frakp \mid d_F\frakD\frakN\}$, which runs in probabilistic polynomial time in the size of the input.  

We now refer to Proposition \ref{conjruntime}
.  We claim that we may arrange so that the diameter $D(G)$ of the $S$-neighbors graph $G=G(S)$ is $D(G)=1$.  The quantity
\[ \frac{\log (M(\frakD,\frakN)w_{\textup{max}}/h^+(\Z_F))}{\log \bigl((N\frakp+1)/(2\sqrt{N\frakp})\bigr)}. \]
has numerator of size $O(\log H)$, so by taking the smallest prime $\frakp \in S$ to be sufficiently large and of polynomial size, we can assume this quantity is $\leq 1$.  

With the set $S$ now computed, we employ Algorithm \ref{Sneighbors}
.  Step 1 can be performed in probabilistic polynomial time over a fixed totally real field $F$ and we have already performed Step $2$.  We conclude from Theorem \ref{conjruntime} that Step 3 requires a number of calls to Problem (\textsf{IsPrincipal}) which is polynomial in the size of the output $H=\#\Cl\calO$, each with input of polynomial size.  In solving Problems (\textsf{IsIsomorphic}) and consequently (\textsf{IsPrincipal}), we only need to check ideals with the same reduced norm, since $S$ consists of primes representing each element in the narrow class group and $D(G)=1$, if $[I] = [J]$ then $[\nrd(I)]=[\nrd(J)]$ so $\nrd(I)=\nrd(J)$.  By Proposition \ref{dembdonn}
, each call to Algorithm \ref{findprincdef} runs in deterministic polynomial time, and the proof is complete.
\end{proof}

We conclude by proposing an alternative to Algorithm \ref{Sneighbors} which solves Problems \ref{computeidealclasses} (\textsf{ClassSet}) and \ref{conjugacyclasses} (\textsf{ConjClassSet}) simultaneously by computing two-sided ideals and connecting ideals, as in Proposition \ref{conjareequiv}.  This algorithm appears to perform better than Algorithm \ref{Sneighbors} in practice, though we do not prove anything rigorous along these lines.

\begin{alg} \label{combined}
Let $\calO$ be an Eichler order of level $\frakN$.  This algorithm outputs a set $\calE=\{\calO_i\}$ of representatives of the conjugacy classes of Eichler orders of level $\frakN$, a set $\calC=\{C_i\}$ of integral invertible right $\calO$-ideals such that $\calO_L(C_i) = \calO_i$, and a set $\mathcal{I}$ of representatives of $\Cl \calO$.
{\rm \begin{enumerate}
\item Solve \textup{(\textsf{ClassNumber})} as in Theorem \ref{existstotdef} and let $H=\#\Cl \calO$.
\item Initialize $\calE, \calC := \{\calO\}$.  Let $\calI$ be a system of representatives for the two-sided $\calO$-ideal classes as in Proposition \ref{alltwosidedclasses}.
\item Choose $\calO_i \in \calE$ and $\frakp$ an ideal of $\Z_F$ coprime to $\frakD\frakN$, and compute the set $\calJ$ of primitive right $\calO_i$-ideals with norm $\frakp$ as in Lemma \ref{idealsnorma}.
\item For all $I \in \calJ$ such that $\calO_L(I)$ is not isomorphic to any order in $\calE$, append $\calO_L(I)$ to $\calE$, and append $\{JIC_i\}_J$ to $\calI$ where $J$ ranges over a set of representatives for the two-sided $\calO_L(I)$-ideal classes.  
\item If $\#\calI = H$, return $\calE,\calC,\calI$; otherwise, return to Step 3.
\end{enumerate}}
\end{alg}

\begin{proof}[Proof of correctness]
The completeness of $\calI$ follows from Proposition \ref{conjareequiv} once we show that the algorithm eventually enumerates all conjugacy classes of orders.  Indeed, let $\calO'$ be an Eichler order of level $\frakN$; then there exists an integral, invertible right $\calO',\calO$-ideal $I$.  By Strong Approximation (Proposition \ref{strongapprox}), we may assume that $\nrd(I)$ is coprime to $\frakD\frakN$.  As in the proof of Lemma \ref{factor_twosided}, $I$ factors into a product $I = I_1 I_2 \cdots I_r$ where each $I_j$ is an invertible integral ideal of prime reduced norm with $\calO_R(I_j)=\calO_L(I_{j+1})$.  The result now follows.
\end{proof}

\begin{rmk}
{\rm In practice, the ability to choose prime ideals $\frakp$ of small norm in Step 3 is essential to enumerate only a small number of ideals of norm $\frakp$; for this reason, we find that it is usually faster to check different orders $\calO_i$ than to compute with only one fixed order as in Algorithm \ref{Sneighbors}.  Nevertheless, the primes $\frakp$ in Step 3 must be chosen in a way such that they generate the narrow class group of $\Z_F$.

Moreover, if some of the orders $\calO_i$ have more than one isomorphism class of two-sided ideals, they contribute to $\mathcal{I}$ accordingly, which speeds up the enumeration.  In this case, the alternative evaluation of the mass formula (Remark \ref{usealtmassformula}) also simplifies, since all right $\calO$-ideals with conjugate left orders have the same mass.  On the other hand, if each $\calO_i$ has only the trivial two-sided ideal class, a condition which is not obviously anticipated, then Algorithm \ref{combined} may take somewhat longer than Algorithm \ref{Sneighbors}.}
\end{rmk}

% \item Let 
% $$b = \min\left\lbrace k \in \mathbb{N} \mid \begin{array}{l} \mbox{for all $1 \le i\le T$ exists an integral right $\calO$-ideal $I_i$}\\ \mbox{with $N(\nrd(I)) \le k$ such that $\calO_i$ is conjugate to $\calO_L(I_i)$}\end{array} \right\rbrace$$
% which does not depend on the choise of the representants $\calO_i$. Then it suffices to check only prime ideals $\frakp$ with $N(\frakp)\le b$.

\section{Definite Eichler orders with class number at most two}

In this section, we list all definite Eichler orders $\calO$ with $h(\calO) \leq 2$.  From (\ref{massdef}) and Proposition \ref{massformula}, for such an order we have
\[ 2 \geq h(\calO) \geq M(\calO) \geq 2^{1-n}|\zeta_F(-1)| \geq \frac{2}{(4\pi^2)^n} d_F^{3/2} \]
and hence 
\begin{equation} \label{rootdiscbound}
\delta_F=d_F^{1/n} \leq (2\pi)^{4/3} \leq 11.594.  
\end{equation}
By the Odlyzko bounds, there are only finitely many such fields and they have been explicitly enumerated \cite{VoightNF}: we have $1,39,47,108,37,40,4,3$ fields of degrees $n=1,2,3,4,5,6,7,8$, respectively, and no field satisfying the bound (\ref{rootdiscbound}) with $n \geq 9$, for a total of $279$ fields.  

For each such field, using the methods of \S 5 (Proposition \ref{computezeta}) we compute the mass $M(\Z_F,\Z_F)$; then for an Eichler order of level $\frakN$ in a quaternion algebra of discriminant $\frakD$ over $F$, we have
\[ M(\frakD,\frakN)=M(\Z_F,\Z_F)\Phi(\frakD)\Psi(\frakN) \leq 2 \]
which gives a finite list of possible ideals $\frakD,\frakN$.  For each such possibility, we compute the corresponding class number using the Eichler mass formula as in Theorem \ref{existstotdef}.  To check the computation, we also enumerate the ideal classes explicitly as in \S 6 (using Algorithm \ref{combined} and the alternate mass formula, Remark \ref{altmassformula}) and see that in all cases they match.  

We consider two Eichler orders, specified by the ideals $\frakD,\frakN$ of $\Z_F$ and $\frakD',\frakN'$ of $\Z_{F'}$, to be \emph{equivalent} if there is an isomorphism $\sigma:F \xrightarrow{\sim} F'$ of fields such that $\sigma(\frakD)=\frakD'$ and $\sigma(\frakN)=\frakN'$.  Two equivalent Eichler orders have the same class number by the mass formula.  

\begin{proposition}
There are exactly $74$ equivalence classes of definite Eichler orders with class number $1$ and $172$ with class number $2$.
\end{proposition}

These classes are listed in Tables \ref{table1}--\ref{table2}.  Here we list the degree $n$, the discriminant $d_F$ of $F$, and the norms $D$ and $N$ of the discriminant $\frakD$ and level $\frakN$.  This way of recording orders is compact but ambiguous; nevertheless, in all cases the field $F$ is determined by its discriminant, and in all but a handful of cases, for \emph{any} choice of squarefree $\frakD$ and coprime $\frakN$ an Eichler order of level $\frakN$ in a quaternion algebra of discriminant $\frakD$ has the given class number.  For the handful of exceptions, we refer to the complete tables which are available online \cite{Voightonline}.

We note that the results in Table \ref{table1} for $F=\Q$ agree with those of Brzezi\'nski \cite{Brzezinski-def} (when restricted to Eichler orders).

\begin{table}[h]
\begin{center}
{\scshape Table \ref{table1}}: Definite Eichler orders $\calO$ with class number $h(\calO)=1$.
\end{center}
\begin{equation} \label{table1} \notag
\begin{array}{cccc|cccc|cccc|cccc}
n & d_F & D & N & n & d_F & D & N &  n  & d_F & D & N & n & d_F & D & N \\
\hline
1 & 1 & 2 & 1 &	2 & 8 & 1 & 1 &	3 & 49 & 7 & 1 &	4 & 725 & 1 & 1 \\
 &  & 2 & 3 &	 &  & 1 & 2 &	 &  & 8 & 1 &	 &  & 1 & 11 \\
 &  & 2 & 5 &	 &  & 1 & 4 &	 &  & 13 & 1 &	 &  & 1 & 19 \\
 &  & 2 & 9 &	 &  & 1 & 7 &	 &  & 29 & 1 &	 &  & 1 & 29 \\
 &  & 2 & 11 &	 &  & 1 & 8 &	 &  & 43 & 1 &	 & 1957 & 1 & 1 \\
 &  & 3 & 1 &	 &  & 1 & 14 &	 & 81 & 3 & 1 &	 &  & 1 & 3 \\
 &  & 3 & 2 &	 &  & 1 & 16 &	 &  & 3 & 8 &	 &  & 1 & 9 \\
 &  & 3 & 4 &	 &  & 1 & 23 &	 &  & 19 & 1 &	 & 2777 & 1 & 1 \\
 &  & 5 & 1 &	 &  & 14 & 1 &	 &  & 37 & 1 &	 &  & 1 & 2 \\
 &  & 5 & 2 &	 &  & 18 & 1 &	 & 148 & 2 & 1 &	 &  & 1 & 4 \\
 &  & 7 & 1 &	 &  & 50 & 1 &	 &  & 2 & 5 &	5 & 24217 & 5 & 1 \\
 &  & 13 & 1 &	 & 13 & 1 & 1 &	 &  & 5 & 1 \\
2 & 5 & 1 & 1 &	 &  & 1 & 3 &	 &  & 5 & 2 &	\\
 &  & 1 & 4 &	 &  & 1 & 9 &	 &  & 13 & 1 &	\\
 &  & 1 & 5 &	 &  & 12 & 1 &	 & 169 & 5 & 1 &	\\
 &  & 1 & 9 &	 & 17 & 1 & 1 &	 &  & 13 & 1 &	\\
 &  & 1 & 11 &	 &  & 1 & 2 &	 & 316 & 2 & 1 &\\
 &  & 1 & 16 &	 &  & 1 & 4 &	 &  & 2 & 2 &\\
 &  & 1 & 19 &	 &  & & 	 &  & 321 & 3 & 1 &\\
 &  & 1 & 20 &		\\
 &  & 1 & 25 &		\\
 &  & 1 & 29 &		\\
 &  & 1 & 44 &		\\
 &  & 1 & 59 &		\\
 &  & 20 & 1 &		\\
 &  & 44 & 1 &		\\
\end{array}
\end{equation}
\end{table}
\addtocounter{equation}{1}

\clearpage

\begin{table}[h]
\begin{center}
{\scshape Table \ref{table2}}: Definite Eichler orders $\calO$ with class number $h(\calO)=2$.
\end{center}
\begin{equation} \label{table2} \notag
\begin{array}{cccc|cccc|cccc|cccc}
n & d_F & D & N & n & d_F & D & N &  n  & d_F & D & N & n & d_F & D & N \\
\hline
1&1&2&7&2&8&1&9&3&49&7&8&4&725&1& 16 \\
&&2&15&   &  &1&17&   &  &7&13&   &  &1& 25 \\
 &  &2&17&&&1&28&   &  &7&27&   &  &1& 31 \\
 &  &2&23&   &  &1&31&   &  &8&7&    &  &1& 41 \\
 &  &3&5&&&1&32&   &  &13&7&   &  &1& 49 \\
 &  &3&7&    &  &1&47&   &  &13&13&  &  &1& 79 \\
 &  &3&8&    &  &14&7&   &  &27&1&   &  &1& 89 \\
 &  &3&11&    &  &34&1&   &  &41&1&   &1125&1& 1 \\
 &  &5&3&   &  &62&1&  &  &71&1&   &  &1& 5 \\
 &  &5&4&    &  &63&1&    &  &97&1&   &  &1& 9 \\
 &  &7&2&    &12&1&1&    &  &113&1&  &  &1& 29 \\
 &  &7&3&    &  &1&2&   &  &127&1&  &  &1& 59 \\
 &  &11&1&&&1&3&   &81&3&17&         &  &80& 1 \\
 &  &17&1&    &  &1&4&    &  &3&19&   &1600&1& 1 \\
 &  &19&1&&&1&6&   &  &8&1&    &1957&1& 7 \\
 &  &30&1&&&1&8&   &  &17&1&   &  &1& 23 \\
 &  &42&1&   &  &1&11&   &  &19&3&   &  &21& 1 \\
 &  &70&1&&&1&12&   &  &73&1&   &2000&20& 1 \\
 &  &78&1&   &  &1&23&  &148&2&17&        &2048&1& 1 \\
2&5&1&31&   &  &6&1&    &  &2&23&   &2225&1& 1 \\
 &  &1&36&   &  &6&11&   &  &5&4&    &2304&18& 1 \\
 &  &1&41&   &  &26&1&   &  &17&1&   &2525&1& 1 \\
 &  &1&45&   &  &39&1&    &  &25&1&   &2624&1& 1 \\
 &  &1&49&   &  &50&1&   &169&5&5&         &2777&1& 8 \\
 &  &1&55&         &13&1&4&   &  &8&1&    &  &1& 11 \\
 &  &1&64&   &  &1&9&  &229&2&1&         &3981&1& 1 \\
 &  &1&71&   &  &1&17&    &  &4&1&    &  &15& 1 \\
 &  &1&79&   &  &1&23&    &  &7&1&    &4205&1& 1 \\
 &  &1&80&   &  &9&1&   &257&3&1&         &4352&14& 1 \\
 &  &1&81&   &  &12&3&   &  &5&1&    &4752&12& 1 \\
 &  &1&89&   &  &39&1&  &  &7&1&    &6809&1& 1 \\
 &  &1&95&   &17&1&4&    &316&2&1&5&14641&11& 1 \\
 &  &1&99&   &  &1&8&    &  &2&4&    &  &23& 1 \\
 &  &20&9&   &  &4&1&   &321&3&1&         &24217&17& 1 \\
 &  &36&1&   &  &18&1&   &  &3&3&    &36497&3& 1 \\
 &  &45&1&   &  &26&1&  &  &7&1&    &38569&7& 1 \\
 &  &55&1&   &21&1&1&   &361&7&1&         &  &13& 1 \\
 &  &95&1&   &  &1&3&  &404&2&1&6&300125&1& 1 \\
 &  &99&1&   &  &1&5&  &469&4&1&         &371293&1& 1 \\
&&124&1&   &  &12&1&  &568&2&1&         &434581&1& 1 \\
&&155&1&   &  &20&1&  &&&&    &485125&1& 1 \\
&&164&1&  &24&6&1& &&&&     &592661&1& 1 \\
&&&&  &  &15&1\\
&&&&  &28&6&1\\
&&&&     &29&1&1\\
&&&&     &33&6&1\\
&&&&     &37&1&1\\
&&&&     &41&1&1
\end{array}
\end{equation}
\end{table}
\addtocounter{equation}{1}

\end{document}